\documentclass{article}

 \usepackage[preprint]{neurips_2026}

\usepackage[T1]{fontenc}    
\usepackage{hyperref}       
\usepackage{url}            
\usepackage{booktabs}       
\usepackage{amsfonts}       
\usepackage{nicefrac}       
\usepackage{microtype}      
\usepackage{xcolor}         
\usepackage{graphicx} 
\usepackage{algorithm}
\usepackage{algorithmic}

\usepackage{natbib}
\usepackage{amsmath}
\usepackage{xspace}
\usepackage{wrapfig}

\definecolor{Pingzhi_color}{HTML}{E0824A}

\definecolor{Shuqing_color}{HTML}{6a994e}

\title{\textbf{\texttt{B$^{3}$-PWL}}: GPU-Batched Branch-and-Bound for Piecewise-Linear Optimization with SOS2 Constraints}

\newcommand{\blfootnote}[1]{%
  \begingroup
  \renewcommand\thefootnote{}\footnote{#1}%
  \addtocounter{footnote}{-1}%
  \endgroup
}

\author{%
  \textbf{Yilin Guan}\,$^{1}$ \quad
  \textbf{Shuqing Luo}\,$^{2}$ \quad
  \textbf{Pingzhi Li}\,$^{2}$ \quad
  \textbf{Tianlong Chen}\,$^{2, \dagger}$ \quad
  \textbf{Kaidi Xu}\,$^{3,4 \dagger}$ \\[0.6em]
  $^{1}$Johns Hopkins University \quad
  $^{2}$University of North Carolina, Chapel Hill \quad \\
  $^{3}$City University of Hong Kong \quad 
  $^{4}$The Hong Kong Institute of AI for Science 
  \\[0.4em]
}

\begin{document}

\maketitle
\blfootnote{$^{\dagger}$Corresponding authors.
  Contact: \texttt{yguan29@jh.edu}, \texttt{tianlong@cs.unc.edu}, \texttt{kaidixu@cityu.edu.hk}.}

\begin{abstract}
  Piecewise-linear (PWL) optimization problems arise in many mixed-integer programming (MIP) optimization applications, including portfolio optimization, workforce scheduling, and resource allocation. But solving them to global optimality remains computationally expensive because branch-and-bound repeatedly solves LP relaxation subproblems. Existing solvers are largely CPU-centric, leaving the scalability of modern GPUs underutilized. Few prior GPU-accelerated branch-and-bound either targets neural network which is not suitable for general PWL optimization, or accelerates only auxiliary subroutines such as strong branching heuristics within CPU-centric MIP solvers. To bridge this gap, we propose \textbf{\texttt{B$^{3}$-PWL}}, a GPU-centric batched branch-and-bound framework for piecewise-linear optimization with Special Ordered Set of type 2 (SOS2) constraints. Our method solves batches of LP relaxation subproblems concurrently on the GPU using a first-order primal-dual solver, enabled by a specialized batched block-tiled sparse matrix kernel. To complement bound computation, we further introduce a unified feasibility search module that combines an SOS2 repair primal heuristic with a batched feasibility pump to rapidly obtain feasible incumbents and improve pruning efficiency. On a benchmark of 43 PWL-MIP instances, \textbf{\texttt{B$^{3}$-PWL}} achieves a 9.25$\times$ geometric-mean speedup over NVIDIA cuOpt while reaching high-quality feasible incumbents on every tested instance. On a public valve-point unit-commitment benchmark, it further outperforms NVIDIA cuOpt and the open-source CPU solvers SCIP and HiGHS, demonstrating the potential of first-order LP methods as the central engine of GPU-accelerated branch-and-bound.
\end{abstract}

\section{Introduction}

Combinatorial optimization with piecewise-linear (PWL) structure underlies many real-world decision-making tasks, including portfolio optimization~\citep{lobo2007portfolio}, workforce scheduling~\citep{ernst2004staff}, and large-scale resource allocation~\citep{karimi2003capacitated}. In mixed-integer programming (MIP), PWL formulations are commonly used to approximate nonlinear relationships while preserving compatibility with optimization frameworks. However, solving such problems to global optimality remains computationally expensive, because the integer constraints induce a combinatorial search space that grows exponentially with problem size. The standard branch-and-bound (B\&B) method tackles this by recursively partitioning the search space and solving an LP relaxation at every node it explores, accumulating thousands to millions of LP solves on a single instance~\citep{land2009automatic, nemhauser1988integer}. In conventional MIP solvers, these LP relaxations are handled by simplex-based methods such as dual simplex~\citep{dantzig1990origins}, whose basis re-optimization and sequential pivoting are difficult to parallelize. As a result, although modern GPUs have utilized massive parallel computing cores to revolutionize large-scale machine learning, the same accelerated computing potential remains largely unexploited for PWL optimization problems.

GPU-accelerated branch-and-bound has been studied in two adjacent settings.
The first line of work targets neural network verification, where frameworks such as $\alpha\beta$-CROWN~\cite{xu2021fast,wang2021beta} run B\&B almost entirely on the GPU. Their efficiency relies on two properties specific to neural network verification. First, neural networks have a layered structure that supports GPU-friendly bound propagation methods such as LiRPA~\cite{xu2020automatic}, which compute linear lower bounds on the network output layer by propagating linear relaxations through the layers, avoiding LP solves at each node. Second, verification only asks whether the worst-case network output violates a safety property, not what its exact value is. Each LP subproblem therefore does not need to be solved exactly: a loose lower bound is enough to prune or further branch the subproblem. This leads to efficient pruning that drastically reduces the number of nodes explored in search tree. However, PWL optimization is more general in structure and it requires solving each LP subproblem precisely to get the globally minimal objective value. Obtaining tight enough bounds and preventing the search tree from exploding are therefore both substantially harder in PWL optimization.

The second line of work uses GPU-accelerated LP solving as an auxiliary component of MIP solvers. NVIDIA cuOpt ~\cite{cuopt, blin2026batched} apply GPU-batched first-order LP methods to specific subroutines such as full strong branching (FSB)~\citep{achterberg2005branching} and optimization-based bound tightening (OBBT)~\cite{gleixner2017three}. However, the LP relaxation at each node of the main search tree is still solved on the CPU using simplex-based methods such as dual simplex. The GPU here only accelerates expensive auxiliary procedures rather than the core B\&B loop.

In this paper, we propose \textbf{\texttt{B$^{3}$-PWL}}, a GPU-batched branch-and-bound framework for piecewise-linear optimization problems with Special Ordered Set of type 2 (SOS2) constraints. Unlike prior work that relegates GPU LP solving to auxiliary subroutines, \textbf{\texttt{B$^{3}$-PWL}} puts batched LP relaxations of subproblems on the GPU as the central compute primitive of the B\&B loop. To support this paradigm, we extend the first-order method cuPDLPx \cite{lu2025cupdlpx} to concurrently solve multiple node-specific LP relaxations within a single batched dispatch, backed by a custom block-tiled sparse matrix kernel that handles the ragged constraint matrices. 
In addition, we introduce a unified feasibility search module that integrate an SOS2 repair primal heuristic~\citep{keha2006branch} and batched feasibility pump strategy tightly and seamlessly into our solving pipeline. It effectively and rapidly repairs and finds feasible solutions to update the best incumbent and tighten the global upper-bound for efficient search tree pruning.

Our contributions are as follows:
\begin{itemize}
\item We present \textbf{\texttt{B$^{3}$-PWL}}, a GPU-batched branch-and-bound framework for piecewise-linear optimization with SOS2 constraints, in which routine node LP relaxations are batched and solved on the GPU via a modified first-order solver for batched LP problems.
\item We develop a unified feasibility search module combining an SOS2 repair primal heuristic and a batched feasibility pump to rapidly obtain feasible incumbents for efficient pruning.
\item On a 43-instance benchmark, \textbf{\texttt{B$^{3}$-PWL}} achieves a 9.25$\times$ geometric-mean speedup and obtain higher quality solutions over cuOpt. We further evaluate on a public valve-point unit-commitment benchmark, where \textbf{\texttt{B$^{3}$-PWL}} outperforms cuOpt, and open-source CPU solvers SCIP and HiGHS, highlighting the potential of first-order LP methods for GPU-accelerated branch-and-bound. 
\end{itemize}

\section{Related Work}

\paragraph{First-order methods for linear programming.}

Simplex algorithm~\citep{dantzig1990origins} and interior-point methods~\citep{wright1997primal, karmarkar1984new} have long dominated linear programming solving. While these algorithms remain the formidable foundation of state-of-the-art CPU solvers, their reliance on matrix factorizations and sequential computations limits their scalability on GPUs. 
To unlock GPU acceleration for large-scale optimization, first-order methods (FOMs) have emerged as a scalable alternative.
PDLP~\citep{applegate2021practical, applegate2025pdlp} reformulates LP as a saddle-point problem and solve it via primal-dual hybrid gradient (PDHG) iterations~\citep{chambolle2011first} with a range of practical enhancements, which transforms FOMs into a practical, numerically stable paradigm for large-scale LPs.
To leverage the massive memory bandwidth of modern accelerators, cuPDLP~\citep{lu2025cupdlp, lu2023cupdlp} extends the CPU-based PDLP to GPUs by offloading core sparse matrix computations and introducing hardware-aware heuristics.
cuPDLPx~\citep{lu2025cupdlpx} further accelerates this paradigm by tailoring the algorithm, achieving multi-fold speedups over its predecessor on large-scale benchmarks. 
Nevertheless, these works focus primarily on accelerating individual LP solves, rather than integrating batched LP computation into B\&B.

\paragraph{GPU-accelerated branch-and-bound.}

GPU acceleration of B\&B has been particularly successful in neural network verification, where $\alpha\beta$-CROWN~\citep{xu2021fast, wang2021beta} pioneered this problem by combining branch-and-bound with GPU-accelerated linear bound propagation (e.g., LiRPA)~\citep{xu2020automatic} and tightening loose bounds with batched gradient descent on GPUs.
In traditional mixed-integer programming (MIP), however, the use of GPUs within exact LP-based branch-and-bound remains more limited. Early efforts incorporating GPUs, such as Neural Branching~\citep{nair2020solving} uses a GPU-scalable ADMM-based variant of full strong branching to generate offline imitation-learning targets for branching policies. More recently, \citet{blin2026batched} extended batched primal-dual hybrid gradient (PDHG)~\citep{chambolle2011first} methods to execute online, exact LP batches on GPUs. However, their application is limited to accelerating specific heuristic sub-routines, including full strong branching (FSB)~\citep{achterberg2005branching} and optimization-based bound tightening (OBBT)~\citep{gleixner2017three}. In contrast, we use batched GPU LP solving as the core engine for routine node processing in exact B\&B for PWL optimization with SOS2 constraints.

\section{Preliminaries}
\label{preliminaries}
\subsection{PWL Problem Formulation with SOS2 Constraints}
We consider an optimization problem in variables $x \in \mathbb{R}^n$, where the first $J$ components $\{x_j\}_{j=1}^J$ are linked to PWL functions. Each PWL function $f_j(x_j)$ is defined over the interval $[\ell_j, u_j]$ with $T_j+1$ ordered breakpoints $\{d_{j,0}, d_{j,1}, \ldots, d_{j,T_j}\}$, partitioning the domain into $T_j$ linear segments with corresponding values $f_j(d_{j,i})$. Using the $\lambda$-formulation with Special Ordered Sets of type 2 (SOS2), we formulate the problem as:
\begin{align}
    \min_{x, z, \lambda} \quad & c^\top x + e^\top z \label{eq:obj} \\
    \text{s.t.} \quad & Ax + Gz \leq b, \quad \ell \leq x \leq u, \notag \\
    & \sum_{i=0}^{T_j} \lambda_{j,i} = 1, \quad \forall\, j \in [J], \label{eq:convex-comb} \\
    & x_j = \sum_{i=0}^{T_j} d_{j,i}\, \lambda_{j,i}, \quad 
      z_j = \sum_{i=0}^{T_j} f_j(d_{j,i})\, \lambda_{j,i}, \quad \forall\, j \in [J], \label{eq:pwl-link} \\
    & \lambda_{j,i} \geq 0, \quad \forall\, i,j, \notag \\
    & \text{SOS2}(\lambda_{j,\cdot}), \quad \forall\, j \in [J]. \notag
\end{align}
Here, $z_j \in \mathbb{R}$ denotes the value of $f_j(x_j)$, and $\lambda_{j,\cdot} = (\lambda_{j,0}, \dots, \lambda_{j,T_j}) \in \mathbb{R}^{T_j+1}_+$ are convex combination weights. 
The SOS2 constraint requires that $\lambda_{j,\cdot}$ to have at most two nonzero components at adjacent indices. That is, if $\lambda_{j,i} > 0$ and $\lambda_{j,i'} > 0$ with $i \neq i'$, then $|i-i'|=1$. This ensures that $(x_j, z_j)$ lies on a single linear segment or exactly on a breakpoint of the PWL graph of $f_j$.

\subsection{Branch and Bound (B\&B) for SOS2 Constraints}
\label{pre:bab}
The combinatorial difficulty of the PWL optimization arises from the SOS2 constraints. Each PWL function $f_j$ admits several candidate active segments, and identifying the correct one for every $f_j$ requires searching over an exponentially large combination of segment choices.
Branch-and-bound (B\&B) solves this problem by recursively restricting the admissible breakpoint range for each PWL function. At the root node, all breakpoints are admissible. Each node $\mathcal{D}$ is characterized by per-function index intervals $[L_j, U_j] \subseteq [0, T_j]$ and adds the restriction
\begin{equation}
    \lambda_{j,i} = 0, \quad \forall\, j \in [J],\; i \notin \{L_j, \ldots, U_j\}.
    \label{eq:node-restriction}
\end{equation}
The subproblem at $\mathcal{D}$ is therefore the original SOS2 formulation \eqref{eq:obj}–\eqref{eq:pwl-link} together with the interval restriction \eqref{eq:node-restriction}. Its LP relaxation, obtained by dropping the SOS2 constraints, provides a lower bound for $\mathcal{D}$. 
Three cases may arise after the LP problem is solved: (i) if the LP is infeasible or its optimum is no better than the incumbent, $\mathcal{D}$ is fathomed; (ii) if the LP solution satisfies all SOS2 constraints, it updates the incumbent if improved and $\mathcal{D}$ is fathomed; (iii) otherwise, $\mathcal{D}$ is branched by splitting one active range into two children. The algorithm terminates when the active queue is empty or the global gap falls within tolerance.

\section{Method}
\label{sec:method}
We introduce \textbf{\texttt{B$^{3}$-PWL}}, a GPU-batched branch-and-bound framework for piecewise-linear optimization with SOS2 constraints. In our framework, batched LP solving is the central compute primitive: the CPU organizes the global search state and assembles task batches, while the GPU executes batches of LP tasks and returns primal and dual solutions. We describe the end-to-end workflow in Section~\ref{sec:overall-framework}, then present the batched node-LP solver in Section~\ref{sec:batched-pdlp}, the block-tiled sparse backend for heterogeneous matrices in Section~\ref{sec:b2-spmv}, and the unified feasibility-search module in Section~\ref{sec:feas-search}.

\subsection{Overall Workflow}
\label{sec:overall-framework}
\begin{figure}[t]
    \centering
    \includegraphics[width=\linewidth]{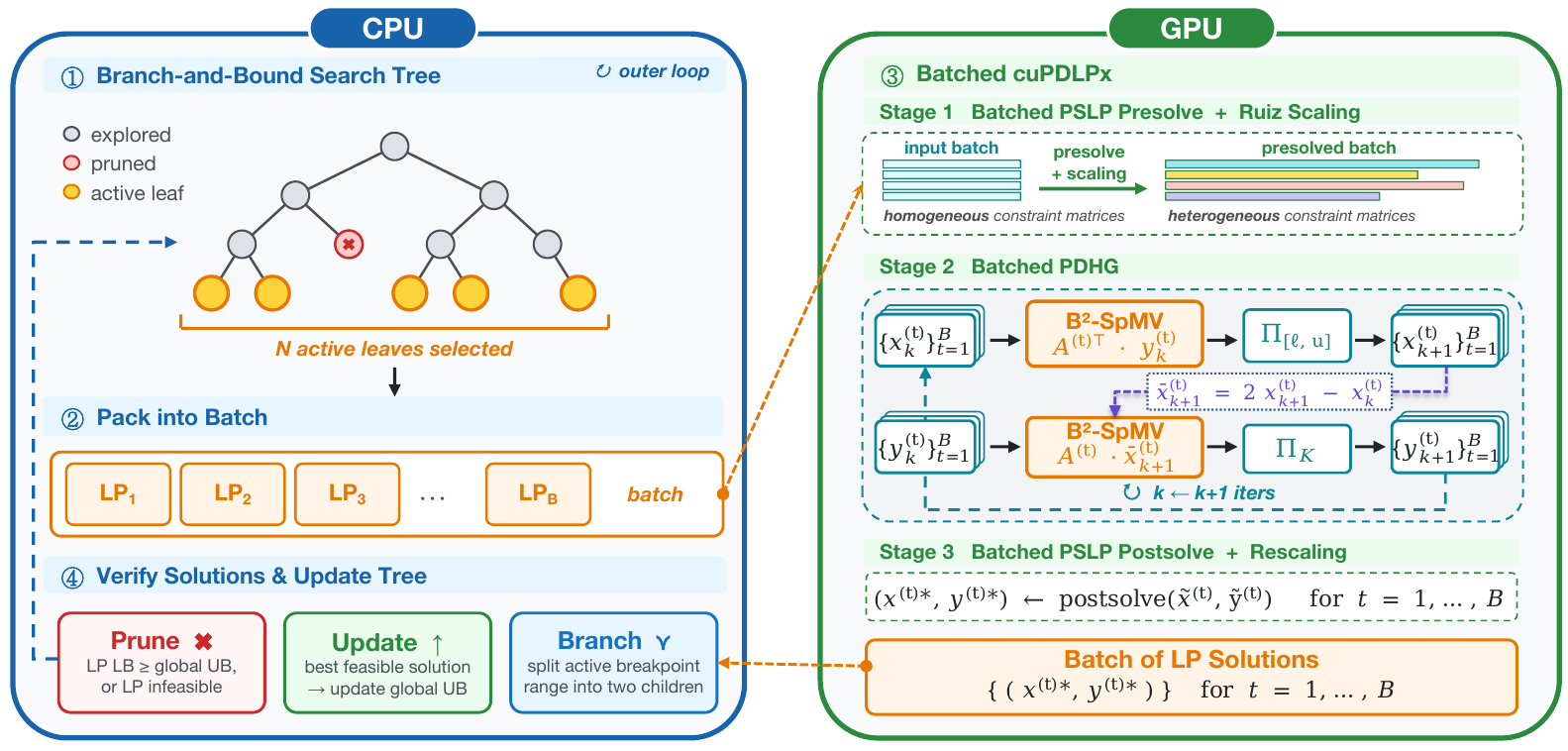}
    \caption{\textbf{Overall workflow for \texttt{B$^{3}$-PWL}.} Each round alternates between CPU-side B\&B management and GPU-side batched LP solving. The CPU (left) selects $N$ active subproblems from the search tree (\textcircled{\small 1}), packs their LP relaxations into a batch (\textcircled{\small 2}), and after the GPU returns solutions, verifies feasibility and updates the tree by pruning, updating the incumbent, or branching (\textcircled{\small 4}). The GPU (right) executes a batched cuPDLPx pipeline (\textcircled{\small 3}): per-task PSLP presolve and Ruiz scaling produce heterogeneous constraint matrices, batched PDHG iterations (built on our B$^2$-SpMV kernel) run until convergence, and postsolve and rescaling yield primal-dual solutions for the entire batch. The outer loop continues until the active queue is empty or the global gap closes.}
    \label{fig:workflow}
\end{figure}

Different node subproblems within a single PWL B\&B tree descend from a common root formulation and differ only in their variable bounds (the active breakpoint ranges), which is naturally suited to batching. Exploiting this, we extend the GPU LP solver cuPDLPx~\citep{lu2025cupdlpx} to concurrently solve multiple node LP relaxations within a single batched dispatch.

Building on this batched solver, our framework proceeds in rounds of three stages. (1) In the \textit{assembly} stage, the CPU selects active nodes from the active nodes queue according to the current B\&B policy and packs their LP relaxations into a single GPU batch. (2) In the \textit{solve} stage, the GPU solves all LP relaxations in the batch in parallel through batched matrix operations, and returns their primal-dual solutions. (3) In the \textit{postprocess} stage, the CPU uses the returned LP solutions to update the global lower and upper bounds, prune fathomed nodes, and branch unresolved nodes into new subproblems that are enqueued into the active node queue. Beyond these standard B\&B operations, we further exploit the LP solutions through a unified feasibility-search module which aims at recovering high-quality feasible incumbents that tighten the upper bound (Section~\ref{sec:feas-search}). The next round then assembles its batch from the updated queue, and the loop continues until the active queue is empty or the global gap falls within tolerance.

\subsection{Heterogeneous Batched LP Engine on GPU}
\label{sec:batched-pdlp}

Each batch consists of $B$ tasks $\mathcal{T}=\{1,\ldots,B\}$, each corresponding to one node LP relaxation. SOS2 constraints are dropped, and the active breakpoint range $[L_j^{(t)},U_j^{(t)}]$ is encoded through the fixings $\lambda_{j,i}^{(t)}=0$ for $i\notin[L_j^{(t)},U_j^{(t)}]$. The objective $c^{(t)}$ is either the original node objective for regular B\&B search tasks or a feasibility-pump objective for feasibility tasks (see Section~\ref{sec:feas-search}). Each task is independently transformed by task-local presolve and Ruiz rescaling into the reduced-space boxed LP form
\begin{equation}
\min_{x^{(t)}}\ c^{(t)\top}x^{(t)}
\quad\text{s.t.}\quad
\ell_A^{(t)} \le A^{(t)}x^{(t)} \le u_A^{(t)},
\qquad
\ell_x^{(t)} \le x^{(t)} \le u_x^{(t)},
\label{eq:reduced_lp}
\end{equation}
yielding reduced matrices $A^{(t)}$ that differ across tasks in dimensions and sparsity patterns. Unlike the prior batched LP engine for heuristic subroutines in cuOpt, where approximate estimation allows presolve to be skipped so that all tasks share a common constraint matrix $A$, our setting requires presolve and therefore motivates a dedicated batched SpMV implementation for heterogeneous matrices $A^{(t)}$ (Section~\ref{sec:b2-spmv}).

Building on cuPDLPx, our engine batches its inner PDHG operator across all $B$ tasks in a single GPU dispatch, while task-local stepsizes, primal weights, projections, restart states, and stopping criteria are maintained independently. The surrounding reflected Halpern updates and adaptive-restart logic follow cuPDLPx and are described in Appendix~\ref{app:cupdlpx}. The inner PDHG operator is dominated by two sparse matrix--vector products.

Let $\mathcal{X}^{(t)}=[\ell_x^{(t)},u_x^{(t)}]$ and $\mathcal{S}^{(t)}=[\ell_A^{(t)},u_A^{(t)}]$, and let $X_k=\{x_k^{(t)}\}_{t=1}^{B}$ and
$Y_k=\{y_k^{(t)}\}_{t=1}^{B}$ denote the batch-wise primal and dual iterates at PDHG iteration $k$. For each task $t$, the inner PDHG operator takes the form
\begin{align}
x_{k+1}^{(t)}
&=
\Pi_{\mathcal{X}^{(t)}}\!\left(
x_k^{(t)}
-
\tau_k^{(t)}
\bigl(
c^{(t)}-A^{(t)\top}y_k^{(t)}
\bigr)
\right),
\label{eq:pdhg-primal}
\\
\bar{x}_{k+1}^{(t)}
&=
2x_{k+1}^{(t)}-x_k^{(t)},
\label{eq:pdhg-extrapolate}
\\
y_{k+1}^{(t)}
&=
y_k^{(t)}
-
\sigma_k^{(t)}A^{(t)}\bar{x}_{k+1}^{(t)}
-
\sigma_k^{(t)}
\Pi_{-\mathcal{S}^{(t)}}\!\left(
(\sigma_k^{(t)})^{-1}y_k^{(t)}
-
A^{(t)}\bar{x}_{k+1}^{(t)}
\right).
\label{eq:pdhg-dual}
\end{align}
Here, $\tau_k^{(t)}=\eta_k^{(t)}/\omega_k^{(t)}$ and $\sigma_k^{(t)}=\eta_k^{(t)}\omega_k^{(t)}$ are the task-local primal and dual stepsizes induced by the cuPDLPx stepsize $\eta_k^{(t)}$ and primal weight $\omega_k^{(t)}$. The row-bound proximal operator projects onto $-\mathcal{S}^{(t)}=[-u_A^{(t)},-\ell_A^{(t)}]$. Full details of cuPDLPx, including the reflected Halpern scheme, adaptive restarting, and primal-weight updates, are provided in Appendix~\ref{app:cupdlpx}.

\subsection{Accelerating Ragged SpMV with Batched Blocking}\label{sec:b2-spmv}

\begin{wrapfigure}{r}{0.45\linewidth}
    \centering
    \vspace{-0.65cm}
    \includegraphics[width=\linewidth]{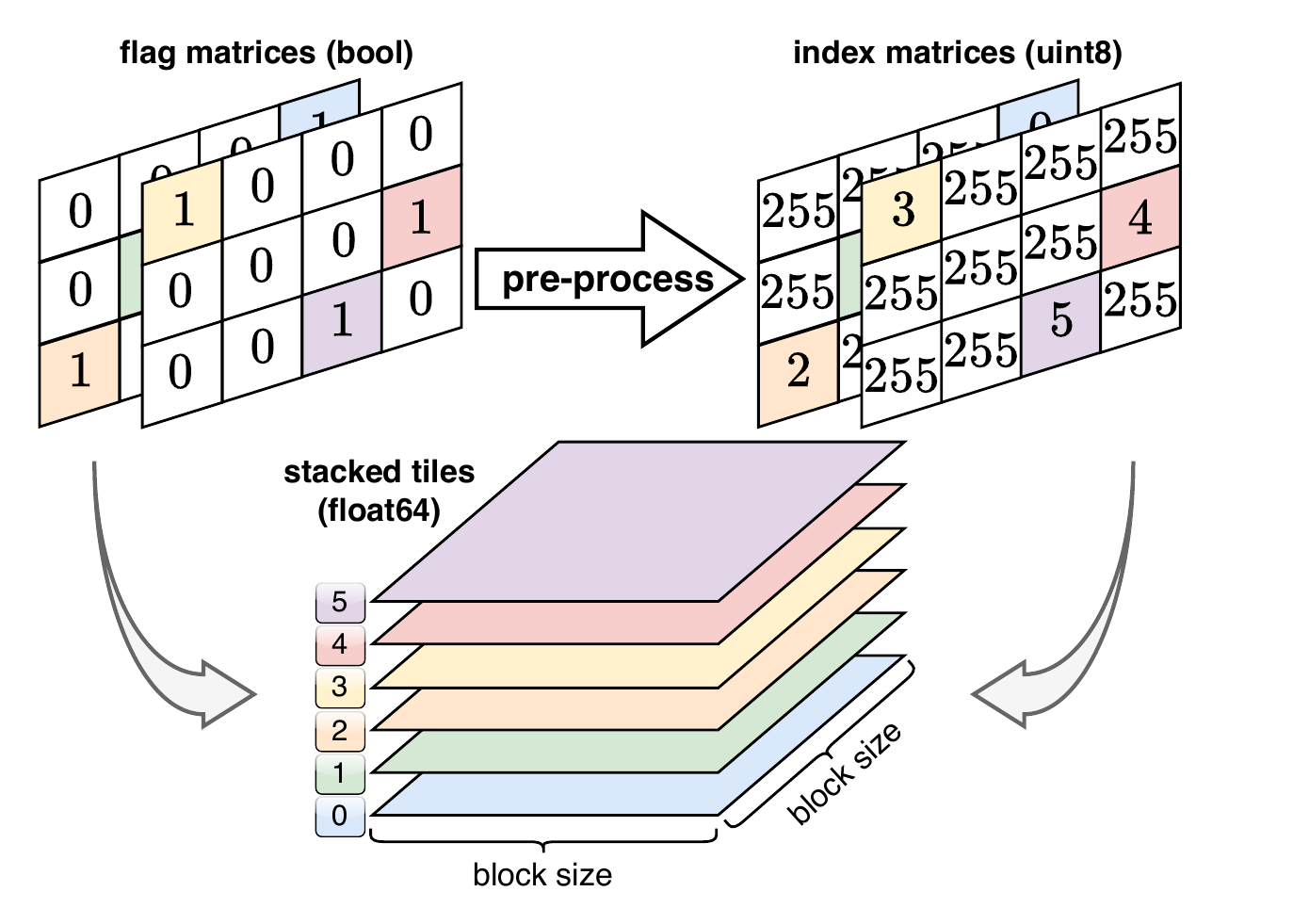}
    \vspace{-0.75cm}
    \caption{\textbf{Building stacked tiles for \texttt{B$^{2}$-SpMV}.}}
    \label{fig:spmv}
    \vspace{-0.5cm}
\end{wrapfigure}

A key challenge in our approach is that different B\&B nodes may produce LP relaxations with different numbers of active variables and constraints after presolve, delivering a \emph{ragged} collection of sparse matrices $\{A^{(t)}\}_{t=1}^B$ with heterogeneous dimensions. Standard GPU-accelerated sparse linear algebra libraries such as cuSPARSE~\cite{cusparse} provide no batched SpMV primitives. Sequentially launching an SpMV kernel per subproblem fails to fully utilize the GPU parallel resources, while flattening the batch into a single block-diagonal sparse matrix leads to poor load balance and inefficient indexing overhead. To support efficient batching on such heterogeneous workloads, we introduce a \emph{batched blocking} scheme for efficient and parallel solving, dubbed \textbf{\texttt{B$^{2}$-SpMV}}. It partitions each sparse matrix $A^{(t)}$ into fixed-size tiles, removes all-zero tiles, and packs the remaining nonzero tiles from all subproblems into a single contiguous GPU tensor. For each subproblem, we build a tile index map that identifies, for every tile position, either the location of the corresponding nonzero tile in the packed tensor or an empty marker, enabling both sparse computing for each LP sub-problem and coalesced memory access for efficient GPU execution.
\texttt{B$^{2}$-SpMV} executes both $A^{(t)} x^{(t)}$ and $A^{(t)\top} y^{(t)}$ with a single batched kernel over the entire ragged batch, where each thread block processes one tiled row or column for one subproblem, with empty tiles skipped based on the index map.

\subsection{Unified Feasibility Search for Efficient Upper-Bound Tightening}\label{sec:feas-search}
Efficient pruning in B\&B depends on tightening either the lower bounds for each node or the global upper bound from the best feasible incumbent. Therefore, high-quality incumbents are critical for pruning by bounding and help drastically reducing the search space. 
However, LP solutions typically violate SOS2 constraints. We therefore introduce a feasibility-search module that attempts to convert LP primal solutions into feasible incumbents through two-stage procedure: a fast and direct SOS2 repair primal heuristic, followed by a batched feasibility pump when direct repair fails.

\paragraph{SOS2 repair primal heuristic.}
Consider an LP solution $(\hat{x}, \hat{z}, \hat{\lambda})$ at node $\mathcal{D}$. We attempt to repair it into an SOS2-feasible one by re-interpolating $\hat{\lambda}_{j,\cdot}$ for each PWL function. The idea is to use $\hat{x}_j$ to identify the active segment it lies on and reset $\hat{\lambda}_{j,\cdot}$ to the convex-combination weights of that segment's two endpoints. 

Concretely, for each PWL function $j$ with active breakpoint range $[L_j, U_j]$, we clip $\hat{x}_j$ to $[d_j^{L_j}, d_j^{U_j}]$, locate the current segment $k_j \in \{L_j, \dots, U_j-1\}$ with $d_j^{k_j} \le \hat{x}_j \le d_j^{k_j+1}$, and reset $\tilde{\lambda}_{j,\cdot}$ as the convex-combination weights of $d_j^{k_j}$ and $d_j^{k_j+1}$:
\[
    \tilde{\lambda}_j^{k_j} = \frac{d_j^{k_j+1} - \hat{x}_j}{d_j^{k_j+1} - d_j^{k_j}}, \qquad \tilde{\lambda}_j^{k_j+1} = \frac{\hat{x}_j - d_j^{k_j}}{d_j^{k_j+1} - d_j^{k_j}},
\]
with $\tilde{\lambda}_j^k = 0$ for all other $k \neq k_j, k_j+1$. The repaired outputs follow as $\tilde{x}_j = \sum_k d_j^k \tilde{\lambda}_j^k$ and $\tilde{z}_j = \sum_k f_j(d_j^k)\, \tilde{\lambda}_j^k$, yielding a candidate $(\tilde{x}, \tilde{z}, \tilde{\lambda})$ that satisfies SOS2 by construction.

Since the repaired candidate satisfies SOS2 by construction, it remains to verify the original linear constraints $A\tilde{x} + G\tilde{z} \le b$ and variable bounds. If satisfied, $(\tilde{x}, \tilde{z}, \tilde{\lambda})$ is feasible for the original PWL problem and updates the incumbent upper bound.

\paragraph{Batched Feasibility-Pump (FP).}

When first-stage repair fails, the repaired vector $\bar{\lambda}^{(r)}$ still encodes which SOS2 segment each $\hat{x}_j$ most naturally suggests. The idea of objective feasibility pump (OFP)~\cite{achterberg2007improving} is to treat $\bar{\lambda}^{(r)}$ as a \emph{soft target} and re-solve the LP under a combined objective that balances the original PWL minimization objective with proximity to $\bar{\lambda}^{(r)}$, so that successive LP solves are pulled toward an SOS2-feasible point while remaining close to optimality, until repair eventually succeeds.

We adapt OFP to our pipeline in two ways: the objective is tailored to the SOS2 formulation, and FP tasks are integrated into the same GPU batch with our B\&B search LP tasks. Specifically, each pump round solves the same node LP with the original constraints intact but a modified objective:
\[
c_{\mathrm{FP}}^{(r)} = \alpha_r \, c_{\mathrm{prox}}^{(r)} + (1-\alpha_r)\, \gamma_r\, c_{\mathrm{orig}}.
\]
where $c_{\mathrm{prox}}^{(r)}$ has entries $(c_{\mathrm{prox}}^{(r)})_{j,k} = -\bar{\lambda}_{j,k}^{(r)}$ on the $\lambda$ block and zero elsewhere, so that minimizing it rewards $\lambda$-mass on the two breakpoints selected by repair, thus pulls the LP solution toward the SOS2 target $\bar{\lambda}^{(r)}$. $c_{\mathrm{orig}}$ is the original objective matrix and keeps the solution close to optimal. $\alpha_r \in [0,1]$ trades off feasibility against quality, and the rescaling factor $\gamma_r = \|c_{\mathrm{prox}}^{(r)}\|_2 / (\|c_{\mathrm{orig}}\|_2 + \varepsilon)$ aligns their $\ell_2$ magnitudes so $\alpha_r$ remains effective. We set $\alpha_0 = 1$ and decay $\alpha_{r+1} = 0.9\, \alpha_r$ on each failure.

Since FP tasks differ from search tasks only in their objective vectors and share the same constraint structure, our batched cuPDLPx solver, which supports per-task objectives, can dispatch them in the same GPU batch as regular B\&B search tasks at no additional scheduling cost. After each FP solve, the returned primal is repaired and checked: accepted candidates update the incumbent, otherwise the repaired point becomes $\bar{\lambda}^{(r+1)}$ for the next round.
\section{Experiments}
\subsection{Experimental Setup}
\label{exp:setup}

\paragraph{Benchmarks.}
Public PWL-MIP instances that keep SOS2 constraints in native form
are limited (detailed in Appendix~\ref{app:benchmarks}), \emph{e.g.} MIPLIB 2017, the most widely used MIP benchmark library, contains no instance with SOS or piecewise-linear constraints, even though the MPS format supports them~\citep{gleixner2021miplib}. Recent PWL work also relies on its own instance generators~\citep{hubner2026spatial}. We therefore evaluate \textbf{\texttt{B$^{3}$-PWL}} on two benchmarks. 

The first is a suite of 43 synthetic PWL-MIP instances from three application families: portfolio optimization, workforce scheduling, and resource allocation. The suite contains 18 portfolio, 9 workforce, and 16 resource-allocation instances, with 100--1,200 PWL functions and 15--25 breakpoints per function. We can control LP size, breakpoint count, sparsity, and SOS2 density through generator to characterize the behavior of our GPU B\&B framework across different scales and structures. Appendix~\ref{app:benchmarks} reports the formulations, generator parameters, seeds, and per-instance size statistics. 

The second is the public valve-point unit commitment benchmark of ~\citet{pedroso2014unit} built on real-world cost data. This benchmark considers thermal unit commitment, where the objective is to determine unit commitment and generation schedules that satisfy multi-period demand while minimizing valve-point fuel costs.

\paragraph{Hardware Settings.}
All end-to-end solver comparisons are run on the same Blackwell node, equipped with NVIDIA RTX PRO 6000 Blackwell Server Edition GPUs (96\,GB each) and a dual-socket x86-64 host with 192 physical cores (384 logical threads) and 718\,GB RAM. 

\paragraph{Solver Configurations.} Both solvers are run with a 1\% relative MIP gap and a 300-second wall-clock limit. \textbf{\texttt{B$^{3}$-PWL}} uses a PDLP feasibility/optimality tolerances of $10^{-4}$. The feasibility-pump module is enabled in the main configuration. cuOpt is run with v26.02.00 which combines CPU-side dual simplex for node LPs, PaPILO presolve, MIR / Gomory / knapsack cutting planes, and GPU-accelerated batched LP for FSB and OBBT.

\paragraph{Metrics.}
Our evaluation focuses on two aspects: solution quality and solve speed.

\textbf{\emph{Solution quality.}}
we use Gurobi as an external oracle reference. We run
Gurobi with tight tolerance gap ($=10^{-6}$) and denote its reference objective by
$z^\star_{\mathrm{ref}}$.
We evaluate solvers on (1) \emph{Exact Match} counts instances on which the solver returns the objective value that exactly matches the solution of Gurobi.
(2) \emph{Eff-Opt} counts instances on
which the solver's objective value $z_{\mathrm{solver}}$ satisfying
\[
    \frac{|z_{\mathrm{solver}} - z^\star_{\mathrm{ref}}|}
    {\max\{1, |z^\star_{\mathrm{ref}}|\}}
    \le 10^{-2}.
\]
Thus, Eff-Opt corresponds to finding a feasible incumbent within $1\%$ of the Gurobi reference objective.

\textbf{\emph{Solve time.}}
All solvers terminate upon reaching a 1\% MIP gap, exhausting the active search queue, or hitting the 300\,s wall-clock budget, whichever comes first; total time is the wall-clock time from solver start to termination. For aggregate runtimes, we report \emph{shifted geometric means} with shift $s=1$\,s following standard MIP benchmarking practice~\citep{achterberg2007constraint, achterberg2009scip}:
\vspace{-0.1em}
\[
    \operatorname{GeoTime}(t)
    =
    \exp\!\left(
        \frac{1}{n}\sum_{i=1}^n \log(t_i+s)
    \right)-s.
\]
\vspace{-0.1em}
The shift reduces the influence of trivially easy instances and avoids overweighting sub-second timing differences. 
\emph{Arithmetic means} report the average of wall-clock times.
\emph{Geo-Spd} is the geometric mean of per-instance
speedup ratios:
\[
    \operatorname{GeoSpd}(\mathrm{cuOpt}/\text{B}^3\text{-PWL})
    =
    \exp\!\left(
        \frac{1}{n}\sum_{i=1}^n
        \log\frac{t_i^{\mathrm{cuOpt}}}{t_i^{\text{B}^3\text{-PWL}}}
    \right).
\]
with values $> 1$ indicating \textbf{\texttt{B$^{3}$-PWL}} is faster on average.
We additionally report \emph{Time to Best Incumbent (TTB)}: the wall-clock time at which a solver first reaches the incumbent it ultimately returns, censored at 300\,s for runs that fail to reach Eff-Opt.
\subsection{Main Results}
\label{exp:main}
Table~\ref{tab:main-results} reports the performance comparisons between \textbf{\texttt{B$^{3}$-PWL}} and cuOpt across our PWL benchmark. Detailed per-instance objectives and wall-clock times are provided in Appendix~\ref{app:per-instance}.

\begin{table}[t]
\centering
\caption{%
  Solution quality and solve time of \textbf{\texttt{B$^{3}$-PWL}} vs. cuOpt on 43 PWL-MIP instances.
}
\label{tab:main-results}
\small
\setlength{\tabcolsep}{3.2pt}
\resizebox{\textwidth}{!}{%
\begin{tabular}{l r cc cc cc cc c c}
\toprule
& & \multicolumn{2}{c}{Eff-Opt} & \multicolumn{2}{c}{Exact}
  & \multicolumn{2}{c}{Arith. Time (s)}
  & \multicolumn{2}{c}{Geo. Time (s)}
  & Geo-Spd & \textbf{\texttt{B$^{3}$-PWL}} \\
\cmidrule(lr){3-4} \cmidrule(lr){5-6}
\cmidrule(lr){7-8} \cmidrule(lr){9-10}
Category & $n$ & \textbf{\texttt{B$^{3}$-PWL}} & cuOpt & \textbf{\texttt{B$^{3}$-PWL}} & cuOpt
  & \textbf{\texttt{B$^{3}$-PWL}} & cuOpt & \textbf{\texttt{B$^{3}$-PWL}} & cuOpt
  & (cu/\textbf{\texttt{B$^{3}$-PWL}}) & Wins \\
\midrule
Portfolio       & 18 & \textbf{18} & 10 & \textbf{15} &  0
  &  2.5 &142.2 &  2.1 & 57.0 & \textbf{29.35}$\times$ & \textbf{18/18} \\
Workforce       &  9 & \textbf{9}  &  8 & \textbf{9}  &  1
  &  2.6 & 68.9 &  2.2 & 35.7 & 17.69$\times$ &  \textbf{9/9} \\
Resource Alloc. & 16 & \textbf{16} & 12 & \textbf{13} &  8
  & 82.3 &107.9 & 39.0 & 67.3 &  1.75$\times$ & 11/16 \\
\midrule
\textbf{Total}  & 43 & \textbf{43} & 30 & \textbf{37} &  9
  & 32.2 &114.1 & \textbf{7.0} & 55.0
  & \textbf{9.25$\times$} & \textbf{38/43} \\
\bottomrule
\end{tabular}%
}
\vspace{-0.1cm}
\end{table}

\begin{table}[t]
\centering
\caption{Comparison with established CPU MIP solvers on the same 43 instances.}
\label{tab:cpu-baselines}
\small
\setlength{\tabcolsep}{8pt}
\begin{tabular}{lcc}
\toprule
Solver & Eff-Opt & Shifted Geo. Time (s) \\
\midrule
\textbf{\texttt{B$^{3}$-PWL}} & \textbf{43/43} & 7.0 \\
Gurobi 13.0.2 & \textbf{43/43} & 1.0 \\
SCIP 10.0.1 & 37/43 & 27.1 \\
HiGHS 1.15.1 & 40/43 & 25.2 \\
cuOpt 26.02 & 30/43 & 55.0 \\
\bottomrule
\end{tabular}
\vspace{-0.3cm}
\end{table}

\paragraph{Solution quality.}
\textbf{\texttt{B$^{3}$-PWL}} reaches Eff-Opt on all 43 instances, compared with 30 / 43 for cuOpt, and matches the Gurobi reference objective (Exact Match) on 37 / 43 instances, compared with 9 / 43 for cuOpt. These results indicate that \textbf{\texttt{B$^{3}$-PWL}} consistently returns higher-quality feasible incumbents than cuOpt on our benchmarks.

\paragraph{Solve speed.}
In total wall-clock time, \textbf{\texttt{B$^{3}$-PWL}} achieves a shifted geometric mean of 7.0 seconds, compared with 55.0 seconds for cuOpt, and an arithmetic mean of 32.2 seconds, compared with 114.1 seconds for cuOpt. Under the reported Geo-Spd metric, this corresponds to a 9.25$\times$ aggregate speedup. We say \textbf{\texttt{B$^{3}$-PWL}} \emph{wins} on an instance if it reaches Eff-Opt while cuOpt does not, or if both reach Eff-Opt and \textbf{\texttt{B$^{3}$-PWL}} terminates faster. Under this criterion \textbf{\texttt{B$^{3}$-PWL}} wins on 38 / 43 instances.

\paragraph{Comparison with established CPU solvers.}
Table~\ref{tab:cpu-baselines} places these results in the context of mature CPU MIP solvers. Gurobi remains the fastest solver on this benchmark, reaching Eff-Opt on all 43 instances with a shifted geometric-mean time of 1.0,s, compared with 7.0,s for \textbf{\texttt{B$^{3}$-PWL}}. \textbf{\texttt{B$^{3}$-PWL}} nevertheless reaches Eff-Opt on all 43 instances and achieves substantially lower aggregate runtime than the open-source CPU solvers SCIP and HiGHS, which reach Eff-Opt on 37/43 and 40/43 instances, respectively. These results position \textbf{\texttt{B$^{3}$-PWL}} as a proof of concept that batched first-order LP solving can serve as the core engine of GPU-centric branch-and-bound, rather than as a replacement for mature commercial MIP solvers.

\subsection{Public valve-point UC benchmark}
We further evaluate \textbf{\texttt{B$^{3}$-PWL}} on the public valve-point unit commitment benchmark of \citet{pedroso2014unit}. Since valve-point costs are nonconvex and nonsmooth, the original paper solves the problem through an adaptive piecewise-linear approximation, yielding a sequence of PWL-MIP instances.
We use the public 40-unit \texttt{eld40}/\texttt{ucp40} data with the original horizons $T \in \{1,3,6,12,24\}$. Unlike the adaptive breakpoint-refinement procedure of the original paper, we fix each cost curve to a static breakpoint count (20 unless noted) to obtain deterministic solver-comparison instances. 

As shown in Table~\ref{tab:vp-time}, Gurobi remains fastest and we do not claim otherwise. Against the open-source simplex-based solvers B$^3$-PWL is faster on 6 of 7 instances (geomean 10.3$\times$ over SCIP, 5.1$\times$ over HiGHS). The
exception is $T=1$, too small to amortize GPU setup. The margin widens with horizon length. At $T=24$, SCIP and cuOpt exhaust the budget at every breakpoint setting while \textbf{\texttt{B$^{3}$-PWL}} certifies the gap in under $10$\,s.

\begin{table}[t]
\centering
\caption{Time to a $1\%$ relative MIP gap (s) on the valve-point UC benchmark.
Rows 1-5 sweep the horizon $T$ at 20 breakpoints; rows 5-7 sweep breakpoints at $T=24$. Geo-Spd is the geometric mean of per-instance ratios against \textbf{\texttt{B$^{3}$-PWL}}.}
\label{tab:vp-time}
\small
\setlength{\tabcolsep}{6pt}
\begin{tabular}{l rrrrr}
\toprule
Instance & \textbf{\texttt{B$^{3}$-PWL}} & Gurobi & SCIP & HiGHS & cuOpt \\
\midrule
$T=1$           & 1.13 & 4.02 &   0.12 &   0.06 &   3.33 \\
$T=3$           & 3.47 & 0.45 &  37.8  & 211.4  & 245.9  \\
$T=6$           & 2.12 & 0.61 &  21.1  &  11.0  &  10.9  \\
$T=12$          & 2.92 & 1.19 &  36.4  &  14.2  & 209.0  \\
$T=24$, bp$=$20 & 5.93 & 3.04 & $>$300 &  39.9  & $>$300 \\
\midrule
$T=24$, bp$=$15 & 5.73 & 2.36 & $>$300 &  34.6  & $>$300 \\
$T=24$, bp$=$30 & 9.02 & 3.95 & $>$300 & 254.0  & $>$300 \\
\midrule
Geo-Spd & --- & 0.47$\times$ & 10.34$\times$ & 5.13$\times$ & 25.40$\times$ \\
\bottomrule
\end{tabular}
\end{table}

\begin{table}[t]
\centering
\caption{Final relative gap (\%) and best objective on the valve-point UC benchmark. ``n/a'': no feasible solution within the budget.}
\label{tab:vp-quality}
\small
\setlength{\tabcolsep}{3.2pt}
\resizebox{\textwidth}{!}{%
\begin{tabular}{l rr rr rr rr rr}
\toprule
& \multicolumn{2}{c}{\textbf{\texttt{B$^{3}$-PWL}}} & \multicolumn{2}{c}{Gurobi}
& \multicolumn{2}{c}{SCIP} & \multicolumn{2}{c}{HiGHS}
& \multicolumn{2}{c}{cuOpt} \\
\cmidrule(lr){2-3}\cmidrule(lr){4-5}\cmidrule(lr){6-7}
\cmidrule(lr){8-9}\cmidrule(lr){10-11}
Instance & Gap & Obj.\ & Gap & Obj.\ & Gap & Obj.\ & Gap & Obj.\ & Gap & Obj.\ \\
\midrule
$T=1$   & 0.052 & 122{,}356 & 0.485 & 122{,}889 & 0.258 & 122{,}609
        & 0.010 & 122{,}305 & 0.086 & 122{,}401 \\
$T=3$   & 0.999 & 293{,}101 & 0.990 & 293{,}375 & 0.929 & 293{,}394
        & 0.460 & 292{,}890 & 0.876 & 293{,}266 \\
$T=6$   & 0.544 & 587{,}487 & 0.371 & 586{,}675 & 0.634 & 588{,}350
        & 0.288 & 586{,}440 & 0.436 & 586{,}981 \\
$T=12$  & 0.327 & 1{,}192{,}044 & 0.673 & 1{,}196{,}353 & 0.827 & 1{,}198{,}228
        & 0.289 & 1{,}191{,}898 & 0.905 & 1{,}199{,}171 \\
$T=24$, bp$=$20 & 0.369 & 2{,}401{,}839 & 0.182 & 2{,}397{,}989
        & 6.781 & 2{,}556{,}152 & 0.262 & 2{,}400{,}451
        & 5.721 & 2{,}538{,}496 \\
\midrule
$T=24$, bp$=$15 & 0.345 & 2{,}404{,}897 & 0.201 & 2{,}402{,}118
        & 6.242 & 2{,}547{,}243 & 0.233 & 2{,}403{,}133
        & 1.682 & 2{,}437{,}913 \\
$T=24$, bp$=$30 & 0.414 & 2{,}398{,}461 & 0.189 & 2{,}393{,}431
        & \multicolumn{2}{c}{n/a} & 0.115 & 2{,}391{,}921
        & 1.894 & 2{,}434{,}827 \\
\bottomrule
\end{tabular}%
}
\vspace{-0.3cm}
\end{table}

\subsection{Time to Best Incumbent}
\label{sec:ttb}

Total wall-clock time conflates two phases of branch-and-bound: finding a high-quality feasible incumbent and closing the final dual-bound gap. We separate them by reporting Time to Best Incumbent (TTB) (defined in Section~\ref{exp:setup}), censored at 300\,s for runs that do not reach Eff-Opt. Since \textbf{\texttt{B$^{3}$-PWL}} reaches Eff-Opt on every instance, its TTB is observed in all 43 cases. cuOpt's TTB is censored on the 13 instances where it fails to reach Eff-Opt within the budget.

Table~\ref{tab:ttb} reports the comparison. Across all 43 instances, \textbf{\texttt{B$^{3}$-PWL}} achieves a $12.28\times$ geometric-mean TTB speedup over cuOpt overall (4.22\,s vs.\ 45.77\,s shifted geomean), and a strictly lower TTB on every instance. Per-category speedups follow the same pattern as total runtime in Table~\ref{tab:main-results}: $28.66\times$ on portfolio and $17.56\times$ on workforce, narrowing to $3.87\times$ on resource allocation, where cuOpt finds incumbents quickly via dual simplex.

Comparing Table~\ref{tab:ttb} with the total solve times in Table~\ref{tab:main-results} reveals \textbf{\texttt{B$^{3}$-PWL}}'s end-to-end advantage stems primarily from rapid primal progress: high-quality incumbents are typically found within the first few GPU batches, after which most of the remaining wall-clock budget is spent on dual-bound tightening to certify optimality. This profile suggests that further end-to-end gains will come from accelerating dual-side convergence such as tighter LP relaxations.

\begin{table}[t]
\centering
\caption{%
Time to best incumbent (TTB) on the 43-instance benchmark. Geo. TTB is the
shifted geometric mean with shift $s=1$\,s; Geo-Spd is the geometric mean of
per-instance cuOpt/\textbf{\texttt{B$^{3}$-PWL}} TTB ratios.
}
\label{tab:ttb}
\small
\setlength{\tabcolsep}{4pt}
\begin{tabular}{l r rr r rr}
\toprule
 & & \multicolumn{2}{c}{Geo.\ TTB (s)} & &
     \multicolumn{2}{c}{Arith.\ TTB (s)} \\
\cmidrule(lr){3-4} \cmidrule(lr){6-7}
Category & $n$ & \textbf{\texttt{B$^{3}$-PWL}} & cuOpt & Geo-Spd & \textbf{\texttt{B$^{3}$-PWL}} & cuOpt \\
\midrule
Portfolio       & 18 & 2.08 &  55.70 & 28.66$\times$ &   2.5 & 141.9 \\
Workforce       &  9 & 2.15 &  35.23 & 17.56$\times$ &   2.5 &  68.4 \\
Resource Alloc. & 16 &11.53 &  42.48 &  3.87$\times$ &  60.4 &  96.0 \\
\midrule
\textbf{Total}  & 43 & \textbf{4.22} & 45.77 & \textbf{12.28$\times$} & \textbf{24.1} & 109.4 \\
\bottomrule
\end{tabular}
\end{table}

\subsection{Unified Feasibility Search Module Ablation}
\label{sec:fp-ablation}
We ablate the two stages of the feasibility-search module in Section~\ref{sec:feas-search}, the SOS2 repair primal heuristic and the batched Feasibility Pump (FP), under three configurations: (i) no feasibility heuristic, (ii) repair only, and (iii) repair + FP. Table~\ref{tab:fp-ablation} reports solve time across the 43-instance benchmark.

\paragraph{SOS2 repair primal heuristic is necessary for PDLP outputs.} Unlike simplex-based solvers, which return vertex solutions that naturally satisfy integrality at leaf nodes, PDHG iterates in PDLP converge from the interior of the LP polytope. Its postsolved primal iterate may satisfy the LP relaxation within global tolerance while still leaving small nonzero mass on SOS2 coefficients that should be exactly zero. Consequently, even leaf nodes whose branching restrictions identify a single admissible SOS2 segment may fail the strict incumbent check. Column~(i) of Table~\ref{tab:fp-ablation} confirms this empirically: no feasible incumbent is recorded on any of the 43 instances. The SOS2 repair step resolves the mismatch by projecting each PWL component onto an admissible segment. With repair enabled (column~(ii)), feasibility is recovered on all 43 instances at low overhead (Table~\ref{tab:time-breakdown}).

\paragraph{FP provides category-dependent acceleration.}
The Feasibility Pump injects proximity-objective LP tasks into the GPU batch to accelerate discovery of high-quality incumbents. Its benefit depends on whether early incumbent quality is the bottleneck for pruning. On portfolio and workforce instances, where the first high-quality incumbent dominates total solve time, FP yields geometric-mean speedups of $2.18\times$ and $2.61\times$ respectively over repair alone. On resource allocation, however, the branching heuristic already discovers incumbents cheaply at nearby leaves, so FP yields no measurable acceleration ($1.00\times$).

\begin{table}[t]
\centering
\caption{Ablation of the feasibility-search module on the 43-instance benchmark. \emph{Feas} reports instances on which a feasible incumbent is recorded; configurations (ii) and (iii) achieve 43/43 on every category. The last column reports the runtime ratio (ii)/(iii); values $> 1$ mean FP accelerates over repair alone.}
\label{tab:fp-ablation}
\small
\setlength{\tabcolsep}{5pt}
\begin{tabular}{lc cc cc cc c}
\toprule
& & \multicolumn{2}{c}{(i) No heuristic} & \multicolumn{2}{c}{(ii) Repair only} & \multicolumn{2}{c}{(iii) Repair + FP} & FP Speedup \\
\cmidrule(lr){3-4} \cmidrule(lr){5-6} \cmidrule(lr){7-8} \cmidrule(lr){9-9}
Category & $n$ & Feas & Geo. (s) & Geo. (s) & Arith. (s) & Geo. (s) & Arith. (s) & (ii)/(iii) \\
\midrule
Portfolio       & 18 & 0/18 & 300.0 & 4.97            & 36.2           & \textbf{2.09}  & \textbf{2.5}  & \textbf{2.18$\times$} \\
Workforce       &  9 & 0/9  & 300.0 & 6.10            & 37.2           & \textbf{2.17}  & \textbf{2.6}  & \textbf{2.61$\times$} \\
Resource Alloc. & 16 & 0/16 &  39.8 & 38.98 & 82.3  & 38.98          & 82.3          & 1.00$\times$ \\
\midrule
Total           & 43 & 0/43 & 142.1 & 11.56 & 53.6  & \textbf{7.05}  & \textbf{32.2} & \textbf{1.69$\times$} \\
\bottomrule
\vspace{-0.5cm}
\end{tabular}
\end{table}

\begin{table}[t]
\centering
\caption{%
Time-weighted breakdown of \textbf{\texttt{B$^{3}$-PWL}}'s branch-and-bound phase. Each entry reports $\sum_i t_i^{\mathrm{component}} / \sum_i t_i^{\mathrm{BaB}}$ within the corresponding category.
}
\label{tab:time-breakdown}
\small
\setlength{\tabcolsep}{4pt}
\begin{tabular}{l rrr r}
\toprule
Component                        & Portfolio & Workforce & Resalloc & Overall \\
\midrule
\textbf{PDHG loop $\&$ GPU setup}  & \textbf{82.1\%} & \textbf{70.4\%} & \textbf{73.9\%} & \textbf{74.0\%} \\
PSLP presolve/rescale            &  9.1\% & 15.9\% & 11.4\% & 11.3\% \\
SOS2 check $\&$ primal heuristic    &  2.0\% &  1.0\% &  6.0\% &  5.9\% \\
Branch selection $\&$ node expansion    &  4.4\% &  4.9\% &  3.7\% &  3.8\% \\
Batch construction               &  3.8\% &  4.7\% &  3.3\% &  3.3\% \\
Other postprocessing             &  0.1\% &  3.2\% &  2.0\% &  2.0\% \\
\bottomrule
\end{tabular}
\end{table}

\subsection{Accelerating Ragged SpMV with Batched Blocking}\label{sec:b2}

\begin{wrapfigure}{r}{0.39\textwidth}
    \centering
    \vspace{-1.2cm}
    \includegraphics[width=0.99\linewidth]{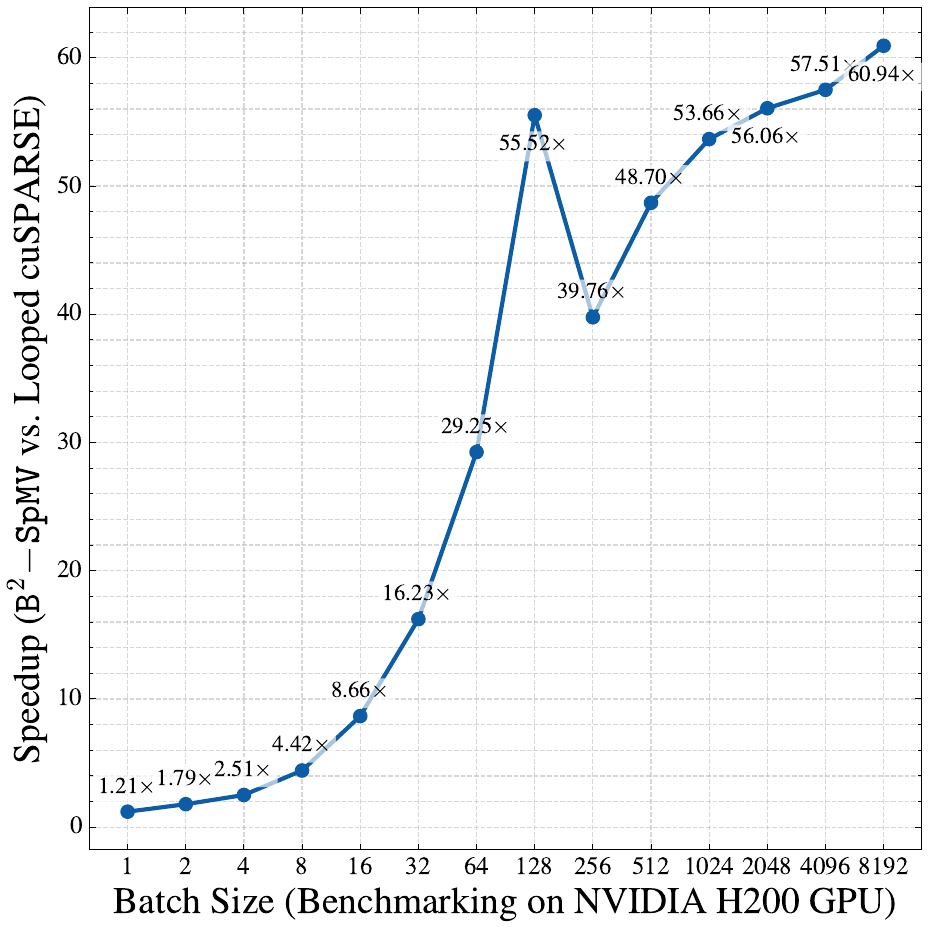}
    \vspace{-0.76cm}
    \caption{\textbf{Speedup Visualization for \texttt{B$^{2}$-SpMV} vs. Looped cuSPARSE}. }
    \label{fig:speedup}
    \vspace{-1.2cm}
\end{wrapfigure}

We benchmark the efficiency of \texttt{B$^{2}$-SpMV} with varied batch size, where the size of the constraint matrix per sample is randomly generated within $8000\sim25000$. We use the default SpMV interface from cuSPARSE as baseline, with looped execution for batching. It scales with batch size and saturates beyond 50$\times$ at larger batches, confirming that the kernel effectively amortizes per-task launch and memory overhead across the batch.

\subsection{Time Breakdown}
\label{sec:time-breakdown}
Table~\ref{tab:time-breakdown} reports the time-weighted breakdown across the 43 instances. Batched LP solving (PDHG loop, task-local PSLP presolve and rescaling) consumes over 85\% of B\&B time, while all CPU-side operations (SOS2 checking, primal heuristics, branching, batch construction, and postprocessing) account for less than 15\% combined. This pattern holds consistently across all three categories, indicating that \textbf{\texttt{B$^{3}$-PWL}} is bound by GPU-side LP throughput rather than CPU-side tree management.

\section{Conclusion}

We introduced \textbf{\texttt{B$^{3}$-PWL}}, a GPU-batched branch-and-bound framework for PWL optimization with SOS2 constraints that solves batched node LP relaxations on the GPU as the main engine of the B\&B loop. On a 43-instance benchmark, \textbf{\texttt{B$^{3}$-PWL}} achieves a 9.25$\times$ geometric-mean speedup over NVIDIA cuOpt while reaching high-quality incumbents on every tested instance. We see this work as a first step toward using batched first-order method for LP as the core engine of B\&B, beyond its current role as auxiliary acceleration. Extending the paradigm to broader PWL and MIP settings, including the incorporation of cutting planes for tighter dual bounds and more general branching heuristics for efficient tree exploration, remains an open direction for future work.

\newpage
\bibliographystyle{abbrvnat}
\bibliography{references}

\newpage
\newpage
\appendix


\section{Broader Impacts}
\label{app:broader-impact}

\texttt{B$^{3}$-PWL} is a foundational systems contribution: a GPU-accelerated framework for solving piecewise-linear mixed-integer optimization problems to global optimality. It is not tied to any specific deployed application, but PWL-MIP formulations underlie a wide range of decision-making tasks, including portfolio optimization, workforce scheduling, energy and resource allocation, and supply-chain planning. Faster and more accessible exact solvers for this class of problems reduce the computational and energy cost of solving real-world planning problems, and lower the barrier for smaller organizations and academic users who lack access to commercial-grade CPU solvers.

We do not foresee direct negative societal impacts from this work. The framework solves classical optimization problems whose objectives and constraints are specified by the user; it does not learn from data, generate content, or make autonomous decisions. As with any general-purpose optimization tool, downstream impact depends on the application, and we view such considerations as properties of mathematical optimization as a whole rather than risks specific to our contribution.

\section{Limitations and Future Work}
\label{app:limitations}

While \texttt{B$^{3}$-PWL} demonstrates that batched first-order LP solvers can serve as the central engine of branch-and-bound for PWL optimization, the current framework has several limitations that point to natural directions for future work.

\paragraph{Slow dual-bound certification.}
First-order methods such as PDHG converge sublinearly in high-accuracy regimes, making dual-bound tightening substantially slower than with simplex-based solvers. Empirically (Section~\ref{sec:ttb}), \texttt{B$^{3}$-PWL} typically locates a high-quality incumbent within the first few GPU batches, but the remaining wall-clock budget is dominated by closing the final dual gap to certify optimality, and on a small number of instances (e.g., the largest resource-allocation cases in Appendix~\ref{app:per-instance}) the 1\% MIP gap is not closed within the 300\,s budget despite a near-optimal incumbent being found early. Hybridizing first-order LP solves with simplex- or interior-point-based crossover for high-accuracy phases is a promising direction for narrowing this gap.

\paragraph{No warm starting between subproblems.}
Simplex-based MIP solvers exploit the structural similarity between parent and child node LPs through dual-simplex warm starting from the parent basis, which often resolves a child node in only a handful of pivots. PDHG, in contrast, lacks a basis representation, and our batched first-order method solvers re-solve each subproblem essentially from scratch despite the high similarity between sibling LPs. Designing warm-start strategies for batched first-order B\&B (e.g., transferring iterate scaling, primal weights, or partial restart history across related subproblems) could substantially reduce per-node work and is, to our knowledge, an open problem in the first-order MIP literature.

\paragraph{Simple branching and no cutting planes.}
PWL optimization is a broad problem class with highly heterogeneous structure across applications. Our current pipeline focuses primarily on GPU-batched node-LP solving and SOS2-aware feasibility search, but it does not yet incorporate the full range of solver machinery used by highly engineered commercial MIP solvers, such as sophisticated branching rules, cut generation, bound tightening, and advanced presolve strategies. Future work could explore branching heuristics that are robust across diverse PWL structures, as well as GPU-compatible cutting-plane and bound-tightening routines that can be integrated into the batched first-order pipeline. By tightening LP relaxations and guiding the search toward more informative regions of the B\&B tree, these extensions may substantially improve pruning efficiency, reduce dual-bound certification time, and narrow the gap between first-order GPU-centric B\&B and mature commercial solvers.

\section{cuPDLPx: GPU-Accelerated First-Order LP Solver}
\label{app:cupdlpx}

Our framework utilizes cuPDLPx~\citep{lu2025cupdlpx} as the linear programming solver of subproblem relaxations in the B\&B search tree, as it is the state-of-the-art, GPU-accelerated first-order solver designed to solve LP problems as an equivalent primal-dual saddle point problem.

\paragraph{Primal-Dual Hybrid Gradient (PDHG)}
The foundational algorithm of cuPDLPx is the Primal-Dual Hybrid Gradient (PDHG). Given primal variables $x$, dual variables $y$, objective vector $c$, and constraint matrix $A$, PDHG iteratively updates the primal-dual iterate
$z^k = [x^k; y^k]$ by alternating between a primal projection step and a dual
update step:
\begin{align}
    x^{k+1} &= \text{proj}_{\mathcal{X}} \Bigl( x^k - \tau(c - A^\top y^k) \Bigr), \label{eq:pdhg_x} \\
    y^{k+1} &= y^k - \sigma A(2x^{k+1} - x^k) - \sigma \text{proj}_{\mathcal{S}} \Bigl( \sigma^{-1}y^k - A(2x^{k+1} - x^k) \Bigr), \label{eq:pdhg_y}
\end{align}
where $\mathcal{X}$ and $\mathcal{S}$ denote the primal and constraint
feasible domains respectively, and $\tau,\sigma$ are the primal and dual step
sizes. cuPDLPx further reparameterizes $\tau, \sigma$ with a base step size $\eta$ and a primal weight $\omega$ (i.e., $\tau = \eta/\omega$, $\sigma = \eta\omega$). Because operations in \eqref{eq:pdhg_x} and \eqref{eq:pdhg_y} are dominated by matrix-vector multiplications, PDHG are well suited for GPU acceleration.

\paragraph{Algorithmic Enhancements in cuPDLPx}
To overcome the slow tail convergence of vanilla PDHG algorithm, cuPDLPx augments this base routine with four advanced mechanisms, achieving robust and accelerated convergence for large-scale LPs.

\textbf{Reflected Halpern Scheme} Instead of standard average schemes, cuPDLPx employs a Reflected Halpern update. It interpolates between an over-relaxed PDHG step and an anchor point $z^0$:
\begin{equation}
    z^{k+1} = \frac{k+1}{k+2} \Bigl( (1+\gamma)\text{PDHG}(z^k) - \gamma z^k \Bigr) + \frac{1}{k+2} z^0,
\end{equation}
where $\gamma \in [0, 1]$ is the reflection parameter. This mechanism enables larger effective step sizes and strictly bounds the fixed-point residual.

\textbf{Adaptive Restart} The anchor $z^0$ is periodically dynamically reset. cuPDLPx monitors the fixed-point error metric $r(z) = \|z - \text{PDHG}(z)\|_P$. When specific decay conditions are met (e.g., sufficient decay or lack of local progress), it restarts by resetting $z^0$ to the current iterate, keeping the search tightly focused on the optimal neighborhood.

\textbf{Constant Stepsize} To maximize parallel efficiency on GPUs, cuPDLPx eschews the sequential adaptive stepsize heuristics typical of CPU solvers. Instead, it adopts a constant stepsize $\eta \approx 1/\|A\|_2$, efficiently approximated via power iteration.

\textbf{PID-Controlled Primal Weight} To balance the progress between primal and dual spaces, the primal weight $\omega$ is dynamically regulated. At each restart occurrence, cuPDLPx updates $\omega$ using a Proportional-Integral-Derivative (PID) controller based on the logarithmic gap between the primal and dual distances to optimality.

\section{Numerical Robustness and Stopping-Criterion Sensitivity}
\label{appendix:numerical-robustness}

We further examine the numerical robustness of \textbf{\texttt{B$^{3}$-PWL}} with respect to the PDLP feasibility/optimality tolerance and the global B\&B stopping criterion. Specifically, we evaluate the feasibility and stability of the returned incumbents, the validity and stability of the final lower bounds, and sensitivity to tighter target MIP gaps.

\paragraph{Incumbent feasibility and stability.}
We rerun all 43 instances with the PDLP feasibility and optimality tolerances tightened from $10^{-4}$ to $10^{-5}$ and $10^{-6}$. Under all three tolerance settings, every instance returns an incumbent within $1\%$ of the Gurobi reference objective. Relative to the $10^{-4}$ baseline, the maximum change in the solver-reported incumbent objective is only $0.0213\%$ under the tighter $10^{-5}$ and $10^{-6}$ tolerances.

We additionally validate each returned incumbent by fixing its integer/SOS2 pattern and re-solving the remaining continuous LP with HiGHS at a $10^{-9}$ feasibility tolerance. Across the three PDLP tolerance settings, this yields $43\times3=129$ incumbent patterns. All 129 patterns yield feasible solutions, with a maximum constraint violation of $2.91\times10^{-11}$. The maximum objective adjustment introduced by this validation is $0.0180\%$ at $10^{-4}$, $0.000457\%$ at $10^{-5}$, and $0.000933\%$ at $10^{-6}$.

\paragraph{Lower-bound stability and validity.}
For each tighter PDLP tolerance, we compare the final lower bounds against the $10^{-4}$ baseline. We also verify
\[
LB \le z^\star_{\mathrm{Gurobi}} \le UB_{\mathrm{validated}}
\]
for every available finite lower bound, where $UB_{\mathrm{validated}}$ denotes the objective of the independently validated incumbent. Table~\ref{tab:pdlp-tolerance-lb} reports instance-wise relative changes in the final lower bounds with respect to the $10^{-4}$ baseline.

\begin{table}[t]
\centering
\caption{Validity and stability of final lower bounds under tighter PDLP feasibility/optimality tolerances. LB shifts are instance-wise relative changes from the $10^{-4}$ baseline.}
\label{tab:pdlp-tolerance-lb}
\small
\setlength{\tabcolsep}{6pt}
\begin{tabular}{lcccc}
\toprule
PDLP tolerance
& Valid finite bounds
& Median LB shift
& 95th-percentile LB shift
& Maximum LB shift \\
\midrule
$10^{-4}$ & 43/43 & Baseline & Baseline & Baseline \\
$10^{-5}$ & 43/43 & $0.00193\%$ & $0.167\%$ & $0.206\%$ \\
$10^{-6}$ & 43/43 & $0.00329\%$ & $0.371\%$ & $0.989\%$ \\
\bottomrule
\end{tabular}
\end{table}

The final lower bounds remain valid for all 43 instances under each tolerance setting and exhibit only small changes relative to the $10^{-4}$ baseline.

\paragraph{Sensitivity to the B\&B stopping criterion.}
We further tighten the target relative MIP gap from $1\%$ to $0.5\%$ and $0.1\%$, while keeping the PDLP tolerance and 300-second wall-clock limit fixed. Table~\ref{tab:mip-gap-sensitivity} reports incumbent and final-LB shifts relative to the $1\%$ baseline.

\begin{table}[t]
\centering
\caption{Sensitivity to the global B\&B target MIP gap. Incumbent and final-LB shifts are measured relative to the $1\%$ configuration.}
\label{tab:mip-gap-sensitivity}
\small
\setlength{\tabcolsep}{5pt}
\resizebox{\textwidth}{!}{%
\begin{tabular}{lccccc}
\toprule
Target MIP gap
& Eff-Opt
& Maximum incumbent shift
& Maximum final-LB shift
& Valid finite bounds
& Time-limit runs \\
\midrule
$1\%$   & 43/43 & Baseline       & Baseline   & 43/43 & 3  \\
$0.5\%$ & 43/43 & $0.000217\%$   & $0.0989\%$ & 43/43 & 7  \\
$0.1\%$ & 43/43 & $0.000217\%$   & $0.102\%$  & 43/43 & 12 \\
\bottomrule
\end{tabular}%
}
\end{table}

Tightening the MIP-gap criterion does not materially change the returned incumbents or final lower bounds. Its primary effect is to increase the cost of dual-bound certification: the number of runs reaching the 300-second limit increases from 3 at a $1\%$ target gap to 12 at a $0.1\%$ target gap. Thus, the main solution-quality conclusions remain stable under stricter stopping criteria, while closing a tighter global optimality gap requires substantially more computation.

\section{Per-Instance Detailed Results}
\label{app:per-instance}

Tables~\ref{tab:per-inst-portfolio}--\ref{tab:per-inst-resalloc}
give per-instance wall-clock time and returned objective for
PWL-BaB and cuOpt on every canonical instance, alongside the
Gurobi certified optimum ($\textsc{MIPGap}{=}10^{-6}$, separate
32-thread Xeon CPU run; Gurobi runtime is omitted because it is
not directly comparable to the GPU-accelerated PWL/cuOpt times).
\textbf{TO} marks the 300\,s budget exhaustion;
\textbf{---} marks no feasible incumbent returned. Each instance is tagged
\texttt{pwl\{N\}\_s\{seed\}}, where \texttt{N} is the number of PWL
functions in the instance and \texttt{seed} is the random-generator
seed. Within resource allocation, two pairs of size groups share the
same \texttt{N} (\texttt{pwl300} and \texttt{pwl500}); we disambiguate
them with the suffixes \texttt{a}/\texttt{b}, whose generator
parameters are listed in Appendix~\ref{app:benchmarks}.

\begin{table}[ht]
\centering
\caption{Per-instance results: Portfolio (18 instances).}
\label{tab:per-inst-portfolio}
\small
\setlength{\tabcolsep}{5pt}
\begin{tabular}{l rr rr r}
\toprule
 & \multicolumn{2}{c}{PWL-BaB} & \multicolumn{2}{c}{cuOpt} & Gurobi \\
\cmidrule(lr){2-3} \cmidrule(lr){4-5}
Instance & Time (s) & Objective & Time (s) & Objective & Objective \\
\midrule
pwl100\_s42              &   0.6 &    -1636.0085 &    TO &           --- &    -1636.0084 \\
pwl100\_s123             &   0.6 &    -1591.3689 &    TO &           --- &    -1591.3686 \\
pwl100\_s999             &   0.5 &    -1702.2300 &    TO &           --- &    -1702.2300 \\
pwl200\_s42              &   1.2 &    -3329.2898 &    TO &    -3103.7171 &    -3329.2898 \\
pwl200\_s123             &   1.7 &    -3493.2999 &    TO &           --- &    -3493.2999 \\
pwl200\_s999             &   1.2 &    -2820.2763 &    TO &    -2752.9085 &    -2820.2763 \\
pwl300\_s42              &   1.4 &    -5132.8143 &   7.8 &    -5092.7688 &    -5132.8144 \\
pwl300\_s123             &   1.3 &    -5312.8877 &   8.8 &    -5270.0022 &    -5313.4625 \\
pwl300\_s999             &   1.3 &    -4326.8367 &  10.1 &    -4314.2420 &    -4326.8434 \\
pwl400\_s42              &   2.6 &    -6729.8250 &  12.9 &    -6717.2791 &    -6729.8211 \\
pwl400\_s123             &   3.0 &    -6930.5113 &  12.0 &    -6908.1034 &    -6930.5117 \\
pwl400\_s999             &   2.1 &    -6186.0892 &  13.8 &    -6150.4538 &    -6189.1684 \\
pwl500\_s42              &   5.8 &    -8793.2572 &  17.8 &    -8752.8616 &    -8793.2572 \\
pwl500\_s123             &   2.9 &    -8844.0370 &  18.7 &    -8822.4380 &    -8844.0369 \\
pwl500\_s999             &   3.2 &    -8396.4138 &    TO &           --- &    -8398.0055 \\
pwl600\_s42              &   3.6 &    -9917.6555 &    TO &           --- &    -9918.4741 \\
pwl600\_s123             &   3.4 &   -10257.0782 &  27.1 &   -10242.9353 &   -10257.0767 \\
pwl600\_s999             &   8.6 &   -10361.3325 &  30.3 &   -10332.3481 &   -10361.3306 \\
\bottomrule
\end{tabular}
\end{table}

\begin{table}[ht]
\centering
\caption{Per-instance results: Workforce (9 instances).}
\label{tab:per-inst-workforce}
\small
\setlength{\tabcolsep}{5pt}
\begin{tabular}{l rr rr r}
\toprule
 & \multicolumn{2}{c}{PWL-BaB} & \multicolumn{2}{c}{cuOpt} & Gurobi \\
\cmidrule(lr){2-3} \cmidrule(lr){4-5}
Instance & Time (s) & Objective & Time (s) & Objective & Objective \\
\midrule
pwl300\_s42          &   0.8 &      281.7799 &  12.5 &      283.5340 &      281.7799 \\
pwl300\_s123         &   0.8 &      279.7645 &  11.6 &      279.7645 &      279.7645 \\
pwl300\_s999         &   0.8 &      277.9249 &  12.2 &      279.3438 &      277.9249 \\
pwl600\_s42          &   2.1 &      563.5598 &  27.9 &      567.8821 &      563.5598 \\
pwl600\_s123         &   1.8 &      559.5290 &  14.8 &      563.2236 &      559.5290 \\
pwl600\_s999         &   2.1 &      555.8498 &  35.8 &      559.9900 &      555.8498 \\
pwl1200\_s42          &   4.8 &     1121.4303 &   107 &     1127.8207 &     1121.4303 \\
pwl1200\_s123         &   5.5 &     1115.9629 &    TO &           --- &     1115.9629 \\
pwl1200\_s999         &   4.2 &     1121.4227 &  98.3 &     1125.9802 &     1121.4227 \\
\bottomrule
\end{tabular}
\end{table}

\begin{table}[ht]
\centering
\caption{Per-instance results: Resource Allocation (16 instances).}
\label{tab:per-inst-resalloc}
\small
\setlength{\tabcolsep}{5pt}
\begin{tabular}{l rr rr r}
\toprule
 & \multicolumn{2}{c}{PWL-BaB} & \multicolumn{2}{c}{cuOpt} & Gurobi \\
\cmidrule(lr){2-3} \cmidrule(lr){4-5}
Instance & Time (s) & Objective & Time (s) & Objective & Objective \\
\midrule
pwl300a\_s42        &  64.1 &      316.2431 &  28.5 &      316.2431 &      316.2431 \\
pwl300a\_s123       &  64.5 &      316.2431 &  38.4 &      316.2431 &      316.2431 \\
pwl300a\_s999       &  64.0 &      316.2431 &  25.4 &      316.2431 &      316.2431 \\
pwl300b\_s123        &  43.4 &      316.2431 &  23.2 &      316.2431 &      316.2431 \\
pwl400\_s42         &   9.0 &      421.6067 &  31.0 &      421.6067 &      421.6067 \\
pwl400\_s123        &   8.9 &      421.6067 &    TO &           --- &      421.6067 \\
pwl400\_s999        &   8.9 &      421.6067 &  31.4 &      425.3778 &      421.6067 \\
pwl500a\_s42       &  40.0 &      527.0718 &  69.0 &      527.0718 &      527.0718 \\
pwl500a\_s123      &  40.0 &      527.0718 &  37.4 &      527.0718 &      527.0718 \\
pwl500a\_s999      &  39.8 &      527.0718 &  81.0 &      527.7716 &      527.0718 \\
pwl500b\_s42        &  11.3 &      527.0083 &  48.2 &      527.0083 &      527.0083 \\
pwl500b\_s123       &  11.4 &      527.0083 &  60.7 &      531.4298 &      527.0083 \\
pwl500b\_s999       &  11.4 &      527.0083 &  52.1 &      531.7667 &      527.0083 \\
pwl1200\_s42         &    TO &     1278.4221 &    TO &           --- &     1278.1789 \\
pwl1200\_s123        &    TO &     1280.2254 &    TO &           --- &     1280.0844 \\
pwl1200\_s999        &    TO &     1287.2700 &    TO &           --- &     1287.0155 \\
\bottomrule
\end{tabular}
\end{table}

\section{Implementation Details}
\label{app:impl}
\paragraph{Multi-level branching for batch filling.}
GPU throughput is maximized when all $B$ slots are occupied, but the B\&B queue may contain fewer than $B$ nodes, particularly in early iterations.
Inspired by the split-depth mechanism in $\alpha\beta$-CROWN~\citep{xu2021fast}, we preemptively branch nodes multiple levels before solving.
Specifically, when the queue size $|\mathcal{Q}| < 0.9B$, we compute a split depth
\begin{equation}
    d_{\text{split}} = \left\lceil \log_2 \frac{0.9B}{\max(1, |\mathcal{Q}|)} \right\rceil,
\end{equation}
and expand each seed node into up to $2^{d_{\text{split}}}$ leaf nodes via BFS-style recursive branching, using the branching heuristic to select variables at each level.
This ensures the GPU batch is consistently saturated.

\paragraph{Adaptive batch sizing.}
We dynamically adjust $B$ based on GPU memory utilization.
After each solving round, we query the available GPU memory and compute a capacity factor $\gamma = M_{\text{target}} / M_{\text{used}}$, where $M_{\text{target}} = 0.85 \cdot M_{\text{total}}$.
The batch size is updated as $B \leftarrow \min(2B,\; \lfloor \gamma B \rfloor)$, subject to a minimum increase threshold, ensuring stable growth without out-of-memory failures.

\section{Benchmark Details}
\label{app:benchmarks}

To our knowledge, there is no broadly adopted, MIPLIB-like benchmark suite specifically for large-scale nonconvex PWL optimization that preserves native SOS2 structure. In particular, the MIPLIB 2017 paper reports that, although its input format permits SOS and piecewise-linear constraints, none of the submitted instances used these structures. Existing PWL studies therefore often rely on generated application-oriented families. Recent external studies have released reproducible generators and instances for nonconvex PWL network-flow and knapsack problems, but these remain generated benchmarks rather than a broadly adopted standard~\citep{hubner2026spatial}. Public application-specific libraries also exist, notably PGLib-UC, which combines convex piecewise-linear production costs with ordinary binary commitment, startup, and shutdown decisions, and therefore targets a different problem scope~\citep{knueven2020mixed}.

Given this landscape, we construct a controlled systems benchmark to characterize the behavior of our GPU B\&B framework, rather than to claim broad industrial representativeness. We selected three classical PWL application classes, portfolio optimization, workforce scheduling, and resource allocation, whose generators allow us to systematically vary LP size, breakpoint count, sparsity, and SOS2 density and thereby characterize the behavior of our GPU B\&B framework across different scales and structures. We remove root-trivial instances, since they primarily benchmark a single LP relaxation rather than the batched B\&B mechanism studied here. Synthetic generation also allows us to control instance difficulty within an informative range, avoiding cases that are either too easy or too difficult to meaningfully distinguish our framework from CPU simplex-based solvers.

Our benchmark suite spans three problem families that exercise complementary structural regimes for branch-and-bound with batched LP relaxations. The portfolio benchmark (Section~\ref{app:portfolio-benchmark}) features dense SOS2 coverage---one SOS2 substructure per asset---coupled through a small number of global linear constraints. The workforce scheduling benchmark (Section~\ref{app:workforce-benchmark}) features a multi-period grid of SOS2 substructures coupled through skill-coverage and budget constraints. The resource allocation benchmark (Section~\ref{app:resource-alloc-benchmark}) concentrates combinatorial structure on a sparse subset of variables under a block-angular linear coupling. All three families share the same univariate PWL cost-curve generator (described in Section~\ref{app:portfolio-benchmark}) so that differences in solver behavior across families reflect coupling structure rather than curve shape.

\subsection{Portfolio Benchmark with Piecewise-Linear Transaction Costs}
\label{app:portfolio-benchmark}

\paragraph{Problem Formulation.}
We consider a portfolio optimization benchmark with separable PWL transaction-cost penalties. This problem models capital allocation across $N$ assets under a global budget, sparse factor-risk limits, and sector exposure caps. The objective balances expected return against a per-asset transaction penalty.

For each asset $j \in \{1,\dots,N\}$, let $w_j \ge 0$ denote the dollar holding and let $z_j$ denote the transaction-cost penalty associated with $w_j$. The objective is
\begin{equation}
\min_{w,z}\; -\sum_{j=1}^N \mu_j w_j + \gamma \sum_{j=1}^N z_j,
\end{equation}
where $\mu_j$ is the expected return of asset $j$, and $\gamma > 0$ controls the trade-off between return and transaction cost.

The portfolio is subject to three classes of linear coupling constraints. First, the total capital invested cannot exceed the budget $C$:
\begin{equation}
\sum_{j=1}^N w_j \le C.
\end{equation}
Second, for each risk factor $i \in \{1,\dots,F\}$, the factor exposure is bounded by a risk budget $B_i^{\mathrm{risk}}$:
\begin{equation}
\sum_{j=1}^N \sigma_{i,j} w_j \le B_i^{\mathrm{risk}},
\end{equation}
where $\sigma_{i,j}$ is the loading of asset $j$ on factor $i$. Third, sector allocations are capped. Let $s(j)$ denote the sector of asset $j$; then for each sector $s \in \{1,\dots,S_{\mathrm{sec}}\}$,
\begin{equation}
\sum_{j:\, s(j)=s} w_j \le B_s^{\mathrm{sec}}.
\end{equation}

Transaction cost is modeled independently for each asset by a univariate PWL function. For asset $j$, let $\{(d_{j,k}, \phi_{j,k})\}_{k=1}^J$ be the breakpoints of its cost curve, where $d_{j,k}$ is a breakpoint on the holding axis and $\phi_{j,k}$ is the corresponding penalty value. We use the standard SOS2 formulation with convex-combination variables $\lambda_{j,k} \in [0,1]$:
\begin{align}
\sum_{k=1}^J \lambda_{j,k} &= 1, \\
w_j &= \sum_{k=1}^J d_{j,k}\lambda_{j,k}, \\
z_j &= \sum_{k=1}^J \phi_{j,k}\lambda_{j,k},
\end{align}
together with an SOS2 constraint on $\lambda_{j,\cdot}$, so that at most two adjacent breakpoints are active. This enforces that $(w_j, z_j)$ lies on the graph of the PWL penalty function.

The complete formulation is
\begin{equation}
\begin{aligned}
\min_{w,z,\lambda}\quad
& -\sum_{j=1}^{N}\mu_j w_j + \gamma \sum_{j=1}^{N} z_j \\
\text{s.t.}\quad
& \sum_{j=1}^{N} w_j \le C, \\
& \sum_{j=1}^{N} \sigma_{i,j} w_j \le B_i^{\mathrm{risk}}, \qquad i=1,\dots,F, \\
& \sum_{j:\,s(j)=s} w_j \le B_s^{\mathrm{sec}}, \qquad s=1,\dots,S_{\mathrm{sec}}, \\
& \sum_{k=1}^{J}\lambda_{j,k}=1, \qquad j=1,\dots,N, \\
& w_j=\sum_{k=1}^{J} d_{j,k}\lambda_{j,k}, \qquad j=1,\dots,N, \\
& z_j=\sum_{k=1}^{J} \phi_{j,k}\lambda_{j,k}, \qquad j=1,\dots,N, \\
& \lambda_{j,\cdot}\ \text{is SOS2}, \qquad j=1,\dots,N, \\
& w_j \in [\underline w_j,\overline w_j], \qquad j=1,\dots,N.
\end{aligned}
\end{equation}

\paragraph{Data Generation.}
The holding interval for each asset is restricted to a tight box around a reference point $m_j$:
\begin{equation}
[\underline w_j,\overline w_j] = [m_j - 0.35\,\mathrm{span}_j,\; m_j + 0.35\,\mathrm{span}_j].
\end{equation}
This narrow box is chosen deliberately so that no single SOS2 segment is admissible at the LP relaxation, ensuring that branch-and-bound must explore multiple segment combinations rather than terminating at the root LP. The breakpoints $(d_{j,k}, \phi_{j,k})$ are sampled from a univariate U-shaped curve generator that combines a low-volume penalty region, a flat economical mid-range, and a steeply rising tail; the same generator is reused across all three benchmark families to isolate the effect of coupling structure.

Expected returns are sampled independently as $\mu_j \sim \mathcal{U}(0.02, 0.15)$, and the global transaction-cost weight is fixed to $\gamma = 0.005$. Asset-wise scales are drawn as $w_j^{\max} \sim \mathcal{U}(50, 500)$. Assets are assigned to sectors in round-robin order, $s(j) = j \bmod S_{\mathrm{sec}}$. For each factor $i$, the loading matrix $\sigma_{i,j}$ has approximately $30\%$ of entries nonzero (each entry independently masked with probability $0.7$), with nonzero values drawn from $\mathcal{U}(0.01, 0.10)$. Let $w_j^{\max,\mathrm{band}} = \overline w_j$ denote the upper end of the tight box. The total budget is set to
\[
C = \alpha_C \sum_j w_j^{\max,\mathrm{band}}, \qquad \alpha_C \sim \mathcal{U}(0.6, 0.8).
\]
Risk budgets are set by
\[
B_i^{\mathrm{risk}} = \max\!\bigl(\alpha_i\,\sigma_i^\top w^{\max,\mathrm{band}},\,1\bigr), \qquad \alpha_i \sim \mathcal{U}(0.8, 1.2),
\]
and sector budgets by
\[
B_s^{\mathrm{sec}} = \alpha_s \sum_{j:\,s(j)=s} w_j^{\max,\mathrm{band}}, \qquad \alpha_s \sim \mathcal{U}(0.75, 0.95).
\]

\paragraph{Problem Size.}
For $N$ assets and $J$ breakpoints per asset, the model contains
\[
N \text{ holding variables } + N \text{ cost variables } + NJ \text{ SOS2 weights } = N(J+2)
\]
variables in total. The number of global coupling constraints is $1 + F + S_{\mathrm{sec}}$, while the PWL representation contributes $3N$ local constraints. This benchmark is therefore characterized by many homogeneous SOS2 substructures coupled through a relatively small set of sparse linear constraints.

\subsection{Workforce Scheduling with Overtime Costs}
\label{app:workforce-benchmark}

\paragraph{Problem Formulation.}
We consider a multi-period workforce scheduling problem with skill coverage and budget constraints. Let $w = 1, \dots, W$ index workers, $p = 1, \dots, P$ index periods, and $t = 1, \dots, T$ index skills. The decision variable $h_{w,p}$ denotes the number of hours assigned to worker $w$ in period $p$, and $z_{w,p}$ denotes the corresponding labor cost. Each worker has a skill coefficient $\kappa_{w,t} \in [0, 1]$ for skill $t$, and each period $p$ requires at least $m_{t,p}$ effective hours of each skill $t$. In addition, each period is subject to a payroll budget $B_p$.

To model overtime and other nonlinear labor costs, each worker--period pair $(w, p)$ is associated with a piecewise-linear cost function $z_{w,p} = f_{w,p}(h_{w,p})$. We represent this function using $J$ breakpoints $\{(d_{w,p,k}, \phi_{w,p,k})\}_{k=1}^J$ and SOS2 variables $\lambda_{w,p,k}$, yielding
\[
\sum_{k=1}^J \lambda_{w,p,k} = 1, \qquad
h_{w,p} = \sum_{k=1}^J d_{w,p,k}\lambda_{w,p,k}, \qquad
z_{w,p} = \sum_{k=1}^J \phi_{w,p,k}\lambda_{w,p,k},
\]
together with an SOS2 constraint on $\lambda_{w,p,\cdot}$.

The full optimization problem is
\[
\begin{aligned}
\min_{h,z,\lambda}\quad & \sum_{w=1}^W\sum_{p=1}^P z_{w,p}\\
\text{s.t.}\quad
& \sum_{w=1}^W \kappa_{w,t} h_{w,p}\ge m_{t,p}, && \forall t,p,\\
& \sum_{w=1}^W r_w h_{w,p}\le B_p, && \forall p,\\
& \sum_{k=1}^J \lambda_{w,p,k}=1, && \forall w,p,\\
& h_{w,p}=\sum_{k=1}^J d_{w,p,k}\lambda_{w,p,k}, && \forall w,p,\\
& z_{w,p}=\sum_{k=1}^J \phi_{w,p,k}\lambda_{w,p,k}, && \forall w,p,\\
& h_{w,p}\in[\underline h_w,\overline h_w],\ \lambda_{w,p,k}\in[0,1], && \forall w,p,k,\\
& \mathrm{SOS2}(\lambda_{w,p,1},\dots,\lambda_{w,p,J}), && \forall w,p.
\end{aligned}
\]

\paragraph{Data Generation.}
Worker-specific maximum hours are sampled as $h_w^{\max} \sim \mathcal{U}(8, 16)$, wage rates as $r_w \sim \mathcal{U}(15, 50)$, and staffing requirements as $m_{t,p} \sim \mathcal{U}(2, 8)$. The skill matrix $\kappa_{w,t}$ has approximately $40\%$ of entries nonzero, with nonzero entries sampled from $\mathcal{U}(0.3, 1.0)$; if a worker draws no nonzero skill, one randomly chosen entry is resampled from $\mathcal{U}(0.3, 1.0)$ to ensure every worker is qualified for at least one skill. Period budgets are set to $0.7$--$0.9$ times a reference full-capacity payroll level. Each worker--period pair receives an independently sampled PWL cost curve from the same generator used in Section~\ref{app:portfolio-benchmark}, with curves associated with the same worker sharing a common hour range determined by that worker's feasible shift length.

\paragraph{Problem Size.}
For $W$ workers, $P$ periods, and $J$ breakpoints per cost curve, the model contains
\[
WP \text{ hour variables } + WP \text{ cost variables } + WPJ \text{ SOS2 weights } = WP(J+2)
\]
variables in total, plus $TP$ skill-coverage constraints, $P$ payroll-budget constraints, and $3WP$ local linear equations from the PWL encoding. In our benchmark suite, the number of local PWL functions ($WP$) typically ranges from $200$ to $2000$.

\subsection{Block-Angular Resource Allocation with Piecewise-Linear Costs}
\label{app:resource-alloc-benchmark}

\paragraph{Problem Formulation.}
We consider a multi-line production-planning benchmark with block-angular resource structure and separable PWL operating-cost penalties on a subset of products. This problem models a manufacturer operating $G$ product groups (e.g., plants or product lines), in which each group draws on its own dedicated equipment (local resources) but also competes for a shared pool of global resources (e.g., raw-material budget, warehouse capacity, logistics throughput). A subset of products carries a non-convex production-cost curve---capturing setup amortization at low volume, economies of scale at medium volume, and overtime or congestion at high volume---which we model as a univariate PWL penalty.

For each group $g \in \{1, \dots, G\}$, let $\mathcal{J}_g$ denote the set of products in that group, with $|\mathcal{J}_g| = K$ and $N = GK$. For each product $j \in \{1, \dots, N\}$ the decision variable $x_j \ge 0$ is the production level. A subset $\mathcal{P} \subseteq \{1, \dots, N\}$ of size $N_{\mathrm{pwl}}$ carries a PWL cost variable $z_j$, and the objective minimizes the total PWL operating cost,
\begin{equation}
\min_{x,z}\; \sum_{j \in \mathcal{P}} z_j.
\end{equation}

The constraint matrix exhibits a block-angular structure: $G$ diagonal local blocks coupled through one shared global block. For each group $g$, $L$ local capacity constraints couple only products in $\mathcal{J}_g$:
\begin{equation}
\sum_{j \in \mathcal{J}_g} a^{\mathrm{loc}}_{g,r,j}\,x_j \le c^{\mathrm{loc}}_{g,r}, \qquad r = 1, \dots, L,
\end{equation}
while $M$ global capacity constraints couple all $N$ products:
\begin{equation}
\sum_{j=1}^{N} a^{\mathrm{glob}}_{m,j}\,x_j \le c^{\mathrm{glob}}_{m}, \qquad m = 1, \dots, M.
\end{equation}
All linear coupling constraints are inequalities. For each $j \in \mathcal{P}$, the PWL operating cost is parameterized by $J$ breakpoints $\{(d_{j,k}, \phi_{j,k})\}_{k=1}^{J}$, encoded with SOS2 convex-combination variables $\lambda_{j,k} \in [0, 1]$:
\begin{equation}
\sum_{k=1}^J \lambda_{j,k} = 1, \qquad
x_j = \sum_{k=1}^J d_{j,k}\lambda_{j,k}, \qquad
z_j = \sum_{k=1}^J \phi_{j,k}\lambda_{j,k},
\end{equation}
together with an SOS2 constraint on $\lambda_{j,\cdot}$.

The complete formulation is
\[
\begin{aligned}
\min_{x,z,\lambda}\quad
& \sum_{j\in\mathcal{P}} z_j \\
\text{s.t.}\quad
& \sum_{j\in\mathcal{J}_g} a^{\mathrm{loc}}_{g,r,j}\,x_j \le c^{\mathrm{loc}}_{g,r}, && g=1,\dots,G,\ r=1,\dots,L, \\
& \sum_{j=1}^{N} a^{\mathrm{glob}}_{m,j}\,x_j \le c^{\mathrm{glob}}_{m}, && m=1,\dots,M, \\
& \sum_{k=1}^J \lambda_{j,k}=1, && j\in\mathcal{P}, \\
& x_j=\sum_{k=1}^J d_{j,k}\lambda_{j,k}, && j\in\mathcal{P}, \\
& z_j=\sum_{k=1}^J \phi_{j,k}\lambda_{j,k}, && j\in\mathcal{P}, \\
& \mathrm{SOS2}(\lambda_{j,1},\dots,\lambda_{j,J}), && j\in\mathcal{P}, \\
& x_j\in[\underline x_j,\overline x_j], && j=1,\dots,N,
\end{aligned}
\]
where $[\underline x_j, \overline x_j]$ is a tight box for $j \in \mathcal{P}$ (specified below) and the wider default range $[0, 100]$ for $j \notin \mathcal{P}$.

\paragraph{Data Generation.}
For each group $g$, the local block $A^{\mathrm{loc}}_g \in \mathbb{R}^{L \times K}$ is generated with target sparsity $\rho_{\mathrm{loc}} = 0.02$; nonzero positions are drawn uniformly at random and nonzero values are sampled independently from $\mathcal{U}(-2, 2)$. The global block $A^{\mathrm{glob}} \in \mathbb{R}^{M \times N}$ is generated with sparsity $\rho_{\mathrm{glob}} = 0.008$ and the same value distribution; allowing both signs in $A^{\mathrm{glob}}$ models global rows that aggregate signed contributions across products (e.g., net demand vs.\ supply of a shared resource). Capacities are sampled as $c^{\mathrm{loc}}_{g,r}, c^{\mathrm{glob}}_m \sim \mathcal{U}(1, 50)$.

The PWL-active subset $\mathcal{P}$ is drawn uniformly at random without replacement from $\{1, \dots, N\}$; breakpoint values $(d_{j,k}, \phi_{j,k})$ are obtained from the same univariate U-shaped curve generator described in Section~\ref{app:portfolio-benchmark}, applied over the nominal production range $[x_{\min}, x_{\max}] = [0.01, 20]$.

To force a non-trivial branch-and-bound search, every $j \in \mathcal{P}$ is restricted to a tight box of width $0.7(x_{\max} - x_{\min})$ around the midpoint of the PWL range,
\[
[\underline x_j, \overline x_j] = [\bar m - 0.35\,\overline{\mathrm{span}},\; \bar m + 0.35\,\overline{\mathrm{span}}],
\]
with $\bar m = (x_{\min} + x_{\max})/2$ and $\overline{\mathrm{span}} = x_{\max} - x_{\min}$; products outside $\mathcal{P}$ take the wider default range $x_j \in [0, 100]$. After the boxes are imposed, the right-hand sides of both blocks are post-processed by adding a uniform slack $s \sim \mathcal{U}(5, 30)$ on top of the worst-case left-hand side $\min_{x \in [\underline x, \overline x]} A_{i,:}\,x$ for each row $i$, so that the linear feasibility region $\{x : Ax \le c, \underline x \le x \le \overline x\}$ is non-empty by construction.

\paragraph{Problem Size.}
For $N = GK$ products, an active subset of size $N_{\mathrm{pwl}}$, and $J$ breakpoints per PWL function, the model has
\[
N \text{ production variables } + N_{\mathrm{pwl}} \text{ cost variables } + N_{\mathrm{pwl}} J \text{ SOS2 weights } = N + N_{\mathrm{pwl}}(J+1)
\]
variables in total. The number of linear coupling constraints is $GL + M$ (block-diagonal local plus dense global), while the PWL representation contributes $3N_{\mathrm{pwl}}$ local linear equations. In contrast to the portfolio and workforce benchmarks---where every primary variable carries an SOS2 substructure---this benchmark concentrates combinatorial structure on only a fraction of the variable set while keeping a sparse, hierarchical linear coupling.

\section{Problem Characteristics and Solver Behavior}
\label{appendix:problem-properties}

We analyze how structural differences across the three problem families affect the behavior of our GPU B\&B framework.

\paragraph{PWL coverage and sparsity.}
Resource-allocation instances embed a relatively small fraction of PWL variables in a much larger sparse block-angular LP. Because branching on one PWL variable modifies only a local part of this structure, its improvement to the global lower bound can be small relative to the cost of solving the full node LP. Consequently, more branching may be required before effective pruning occurs. Meanwhile, this sparse structure favors mature simplex-based solvers through basis reuse, sparse factorization, cuts, and other B\&B enhancements.

\paragraph{LP size.}
The node LPs in the resource-allocation family are also substantially larger, reducing the maximum batch size that fits within the same GPU memory and increasing the cost of each batched iteration.

\paragraph{Search-tree structure.}
These properties tend to produce narrower, deeper search trees because each local PWL branch provides limited global tightening. In contrast, portfolio and workforce expose many independent PWL decisions and tend to generate shallower, wider frontiers that better match GPU batching. Branching quality is consequently more important for resource allocation: an unproductive early decision can send the search deep into a poor direction and require many weak or infeasible node LPs to be solved before pruning. Unlike simplex methods, first-order methods do not return basic vertex solutions and converge slowly at high dual accuracy, weakening lower-bound-based pruning, especially in deeper search trees.


\end{document}